\documentclass[a4paper]{amsart}

\usepackage{amsmath}
\usepackage{amssymb}
\usepackage[all]{xy}

\usepackage{bbm}
\usepackage{mathrsfs}

\theoremstyle{plain}
\newtheorem{thm}{Theorem}[section]
\newtheorem{lem}[thm]{Lemma}
\newtheorem{prop}[thm]{Proposition}
\newtheorem{cor}[thm]{Corollary}
\newtheorem{ques}[thm]{Question}

\theoremstyle{definition}
\newtheorem{definition}[thm]{Definition}
\newtheorem*{Ack}{Acknowledgements}

\theoremstyle{remark}
\newtheorem{remark}[thm]{Remark}

\newcommand{\kxC}{\mathbbm{C}}  %komplexes C
\newcommand{\natN}{\mathbbm{N}} %natuerliche Zahlen
\newcommand{\prP}{\mathbbm{P}}  %projektiver Raum
\newcommand{\ratQ}{\mathbbm{Q}} %rationale Zahlen
\newcommand{\reR}{\mathbbm{R}}  %reelle Zahlen
\newcommand{\intZ}{\mathbbm{Z}} %ganze Zahlen

\newcommand{\xxD}{\mathcal{D}}
\newcommand{\xxF}{\mathcal{F}} 
 
\newcommand{\xxH}{\mathcal{H}}
\newcommand{\III}{\mathcal{I}}
\newcommand{\xMM}{\mathcal{M}}
\newcommand{\ssh}{\mathcal{O}} 
\newcommand{\xxFf}{\mathscr{F}}

\newcommand{\XXX}{\mathcal{X}} 
\newcommand{\xxZ}{\mathcal{Z}} 
\newcommand{\xLl}{\mathscr{L}}

\newcommand{\Coh}{\operatorname{Coh}}

\newcommand{\xxHom}{\operatorname{Hom}}
\newcommand{\xid}{\operatorname{id}}
\newcommand{\imag}{\operatorname{im}}

\newcommand{\Irr}{\operatorname{Irr}}

\newcommand{\Isom}{\operatorname{\mathbbm{I}som}}

\newcommand{\Mor}{\operatorname{Mor}}

\newcommand{\Spec}{\operatorname{Spec}}
\newcommand{\supp}{\operatorname{supp}}

\newcommand{\tensor}{\otimes_{\kxC}}
\newcommand{\red}{/\!/}

\newcommand{\longhookrightarrow}{\ensuremath{\lhook\joinrel\relbar\joinrel\rightarrow}}
\newcommand{\sgeq}{\mathbin{\geq\hspace{-1.06em}_{_(\:\,_)}}}
\newcommand{\sleq}{\mathbin{\leq\hspace{-1.06em}_{_(\:\,_)}}}

\newcommand{\xxGHilb}[1]{G\operatorname{-Hilb}(#1)} %G-Hilbertschema mit allgemeinem G
\newcommand{\SlHilb}[1]{Sl_2\operatorname{-Hilb}(#1)}

\newcommand{\invHilb}[3]{\operatorname{Hilb}^{#1}_{#2}(#3)} %invariantes Hilbertschema

\newcommand{\finvHilb}[1]{\mathcal{H}ilb_{#1}^G}
\newcommand{\xxGrass}[2]{\operatorname{Grass}({#1},{#2})} %Grassmannsche von Unterrￜumen = {\operatorname{Gr}_{#1}(#2)}
\newcommand{\qGrass}[2]{\operatorname{Grass}({#1},{#2})} %Grassmannsche von Quotienten = {\operatorname{Gr}^{#2}(#1)}
\newcommand{\Quot}{\operatorname{Quot}}
\newcommand{\xxGQuot}[2]{\operatorname{Quot}^G(#1,#2)}
\newcommand{\relGQuot}[3]{\operatorname{Quot}^G_{#1}(#2,#3)}
\newcommand{\xxGQuotH}{\operatorname{Quot}^G(\xxH,h)}
\newcommand{\fGQuotH}{\mathcal{Q}uot^G(\xxH,h)}
\newcommand{\fMck}{\overline{\mathcal{M}}_{\chi,\kappa}(X)}
\newcommand{\fsMck}{\mathcal{M}_{\chi,\kappa}(X)}

\newcommand{\fMth}{\overline{\mathcal{M}}_{\theta}(X)}
\newcommand{\fsMth}{\mathcal{M}_{\theta}(X)}
\newcommand{\Mth}{\overline{M}_{\theta}(X)}
\newcommand{\sMth}{M_{\theta}(X)}

\begin{document}

\title[Moduli spaces of $(G,h)$--constellations]{Moduli spaces of $(G,h)$--constellations}

\author{Tanja Becker}
\address{Tanja Becker\\Institut f\"ur Mathematik\\ Johannes Gutenberg-Universit\"at Mainz\\
Staudingerweg 9\\55099 Mainz\\Germany}

\author{Ronan Terpereau}
\address{Ronan Terpereau\\Institut f\"ur Mathematik\\ Johannes Gutenberg-Universit\"at Mainz\\
Staudingerweg 9\\55099 Mainz\\Germany}
\email{rterpere@uni-mainz.de}

%\received{??}
%\accepted{??}

\begin{abstract}
Given an infinite reductive group $G$ acting on an affine scheme $X$ over $\kxC$ and a Hilbert function $h\colon \Irr G \to \natN_0$, we construct the moduli space $\sMth$ of $\theta$--stable $(G,h)$--constellations on $X$, which is a generalization of the invariant Hilbert scheme after Alexeev and Brion \cite{AB:2005} and an analogue of the moduli space of $\theta$--stable $G$--constellations for finite groups $G$ introduced by Craw and Ishii \cite{CI:2004}. Our construction of a morphism $\sMth \to X \red G$ makes this moduli space a candidate for a resolution of singularities of the quotient $X \red G$.
\end{abstract}

\maketitle

% *** 1 ***************************************************************************************************************************************
\section*{Introduction}

In the study of the action of a reductive group $G$ on an affine scheme of finite type $X$, see \cite{Bec:2010, Tera:2014, Terb:2014}, we observed many situations where the quotient $X \red G$ has at least two different resolutions of singularities $Y_1 \to X \red G$, $Y_2 \to X \red G$, which are dominated by a third one $Z \to X \red G$, for instance as a flop:
%\vspace{-3mm}
$$
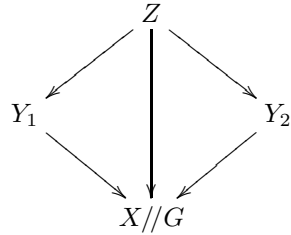
\begin{xy} \xymatrix{
 & Z \ar@{->}[ld] \ar@{->}[dd] \ar@{->}[rd] & \\
Y_1 \ar@{->}[dr] & & Y_2 \ar@{->}[dl]\\
 & X \red G &
}\end{xy} $$

This is for example the case for the action of $Sl_2$ on $(\kxC^2)^{\oplus 6}$ by multiplication from the left and on its subscheme $\mu^{-1}(0)$ defined as the zero fibre of the moment map of the action. Exploring this example further, in \cite{Bec:2011} we found out that in this case the resolution $Z$ is given as an invariant Hilbert scheme $Z = \SlHilb{\mu^{-1}(0)}$. For other similar examples with classical groups acting on $\mu^{-1}(0)$, see \cite{Terb:2014}, and for examples with classical groups acting on classical representations, see \cite{Tera:2014}.   
In \cite{CI:2004}, Craw and Ishii examine this phenomenon for finite groups by introducing the notion of a $G$--constellation on $\kxC^n$ and a certain stability condition $\theta$ on them and by constructing the moduli space $M_{\theta}$ of $\theta$--stable $G$--constellations. For a finite abelian group $G \subset Sl_3(\kxC)$ they show that every projective crepant resolution of $\kxC^3/G$ can be obtained as such a moduli space. Further, choosing $\theta$ appropriately they recover Ito and Nakamura's $G$--Hilbert \linebreak scheme \cite{IN:1996, IN:1999, Nak:2001}.

In this article, we define similar concepts in the case of infinite reductive groups acting on affine schemes of finite type and we construct the moduli space $\sMth$ of $\theta$--stable $(G,h)$--constellations, where the map $h\colon \Irr G \to \natN_0$ replaces the regular representation occurring in \cite{CI:2004}. Let us denote by $\rho_0$ the trivial representation of $G$. If $h(\rho_0)=1$, then for a special choice of $\theta$ we recover Alexeev and Brion's invariant Hilbert scheme \cite{AB:2004, AB:2005, Bri:2010}. Thus, $\sMth$ is a generalization of the invariant Hilbert scheme and an analogue of the moduli space of $G$--constellations, which in turn both generalize the $G$--Hilbert scheme:

Given a finite group $G$ acting on $X$, the $G$--Hilbert scheme parametrizes \linebreak $G$--clusters, i.e.~$G$--invariant subschemes $Z$ of $X$ such that $H^0(\ssh_Z)$ is isomorphic to the regular representation $R$ of $G$. The notion of $G$--constellation introduced in \cite{CI:2004} generalizes this concept when $X=\kxC^n$: a $G$--constellation is a $G$--equivariant coherent $\ssh_X$--module with isotypic decomposition isomorphic to $R$.
Such a $G$--con\-stellation $\xxF$ is $\theta$--stable for some $\theta \in \xxHom_{\natN}(R(G),\ratQ)$, where $R(G)$ is the representation monoid of $G$, if $\theta(\xxF) = 0$ and if for every non--zero proper $G$--equivariant coherent subsheaf $0 \neq \xxF' \subsetneq \xxF$ one has $\theta(\xxF')>0$. In this situation, Craw and Ishii construct the moduli space $M_{\theta}$ of $\theta$--stable $G$--con\-stellations as the GIT--quotient of the space of quiver representations associated to $G$ by the group of $G$--equivariant automorphisms of $R$ as described by King in \cite{Kin:1994}. For a special choice of $\theta$ they recover $M_{\theta} = \xxGHilb{X}$.

A second generalization of the $G$--Hilbert scheme was established by Alexeev and Brion. Fix a complex reductive group $G$, an affine $G$--scheme of finite type $X$, and a map $h\colon \Irr G \to \natN_0$ on the set $\Irr G = \{\rho \colon G \to Gl(V_{\rho})\}$ of isomorphy classes of irreducible representations of $G$. In \cite{AB:2004, AB:2005}, the authors define the invariant Hilbert scheme $\invHilb{G}{h}{X}$, whose closed points parametrize all \linebreak $G$--invariant subschemes of $X$ whose coordinate rings have isotypic decomposition isomorphic to $\bigoplus_{\rho \in \Irr G}\kxC^{h(\rho)}\tensor V_{\rho}$, or equivalently all quotients $\ssh_X/\III$, where $\III$ is an ideal sheaf in $\ssh_X$, with this prescribed isotypic decomposition.

Our contribution to these constructions of moduli spaces is to unify the ideas of \cite{CI:2004} and \cite{AB:2004, AB:2005}: For an infinite complex reductive group $G$, an affine $G$--scheme of finite type $X$ and a map $h\colon \Irr G \to \natN_0$, we define a $(G,h)$--constellation as a $G$--equivariant coherent $\ssh_X$--module with isotypic decomposition given by $h$ as above. Then we introduce $\theta$--stability similarly to the case of $G$--constellations. This stability condition is more delicate than the one of Craw and Ishii since it involves infinitely many parameters. We locate finitely many of them which control the others whenever $\theta$ has finitely many negative values. Then we construct the moduli space of $\theta$--stable $(G,h)$--constellations by means of geometric invariant theory and invariant Quot schemes in a parallel way to the construction of the moduli space of stable vector bundles of Simpson \cite{Sim:1994}.
As a generalization of the Hilbert--Chow morphism we moreover construct a morphism $\sMth \to X \red G$ when $h(\rho_0) = 1$ for $\rho_0$ the trivial representation. Further studies of $\sMth$ have to be made in order to decide in which cases this morphism gives a resolution of singularities.

This article is structured as follows: In Section \ref{chGconst} we set up our framework by introducing the notions of $(G,h)$--constellation, $\theta$--semistability and $\theta$--stability similarly to the case of $G$--constel\-lations and by defining the corresponding moduli functors $\fMth$ and $\fsMth$. Assuming $\theta \in \xxHom_{\natN}(R(G),\ratQ)$ has only finitely many negative values, we show that every $\theta$--stable $(G,h)$--constellation is generated as an $\ssh_X$--module by its components indexed by a certain finite subset $D_- \subset \Irr G$, so that each $\theta$--stable $(G,h)$--constellation is a quotient of a fixed coherent sheaf $\xxH$ and hence an element of the invariant Quot scheme $\xxGQuotH$. With a slightly more restrictive choice of $\theta$, the same holds for $\theta$--semistability.
At the end of this section we show that if $h(\rho_0)=1$ and $\theta_{\rho_0}$ is the only negative value of $\theta$, then the moduli functor $\fsMth$ equals the Hilbert functor $\finvHilb{h}(X)$.

In Section \ref{chQuotGIT} we deal with the geometric invariant theory of the invariant Quot scheme $\xxGQuotH$ in order to construct the moduli space of $(G,h)$--constellations as its GIT--quotient: The invariant Quot scheme is equipped with a certain ample line bundle $\xLl$ coming from an embedding into a product of Grassmannians as established in Subsection \ref{invQuot}. Considering the gauge group $\Gamma$, we examine GIT--stability and GIT--semistability on $\xxGQuotH$ with respect to the induced linearization on $\xLl$ twisted by a certain character $\chi$. Thus, on the set of GIT--semistable quotients $\xxGQuotH^{ss}$ we obtain the categorical quotient $\xxGQuotH^{ss}\red_{\xLl_{\chi}} \Gamma$, which turns out to be a moduli space of GIT--semistable $(G,h)$--constellations in Section \ref{modspace}. 

In Section \ref{chcorres} we establish a correspondence of $(G,h)$--constellations and $G$--equi\-variant quotients $[q\colon \xxH \twoheadrightarrow \xxF] \in \xxGQuotH$ and a correspondence of their respective subobjects. This allows us to introduce another (semi)stability condition $\widetilde \theta$ which is equivalent to GIT--(semi)stability but resembles very much $\theta$--(semi)\-stability. We show that if $\xxF$ is $\theta$--stable, then it is also $\widetilde \theta$--stable and hence any corresponding point $[q\colon \xxH \twoheadrightarrow \xxF]$ in $\xxGQuotH$ is GIT--stable. This allows us to realize the functor $\fsMth$ of flat families of $\theta$--stable $(G,h)$--constellations as a subfunctor of the functor $\fsMck$ of flat families of GIT--stable $(G,h)$--constellations. The same does not work for $\fMth$ and semistability. This shows that the passage from finite to infinite groups is a profound issue.

In Section \ref{modspace} we consider properties of these functors. First we prove that $\theta$--sta\-bility is open in flat families. From this fact we deduce that $\xxGQuotH^{s}_{\theta}$ is an open subscheme of $\xxGQuotH^s$ and hence a quasiprojective scheme. Then we show that the functors $\fMck$ and $\fsMck$ are corepresented by the schemes $\xxGQuotH^{ss}\red_{\xLl_{\chi}} \Gamma$ and $\xxGQuotH^s/\Gamma$, respectively. Denoting by $\xxGQuotH^{s}_{\theta}$ the set of $\theta$--stable elements in $\xxGQuotH$, we also show that $\fsMth$ is corepresented by the scheme $\xxGQuotH^{s}_{\theta}/\Gamma$. We call 
$$\sMth := \xxGQuotH^{s}_{\theta}/\Gamma$$ 
the moduli space of $\theta$--stable $(G,h)$--constellations.  We define the scheme $\Mth$ as the closure of $\sMth$ in $\xxGQuotH^{ss}\red_{\xLl_{\chi}} \Gamma$. Finally, when $h(\rho_0) = 1$ we construct a morphism from $\Mth$ to the quotient $X\red G$ corresponding to the Hilbert--Chow morphism.

In the outlook (Section \ref{secO}) we discuss some further aspects of the moduli spaces $\sMth$ and $\Mth$, which are worth pursuing in the future.

\begin{Ack}
The present article grew up in the first-named author's thesis \cite{Bec:thesis} and was then completed by the second-named author. Some proofs are omitted here since they are analogous to existing literature or they consist of easy calculations worked out in detail in \cite{Bec:thesis}. 
The first-named author thanks her advisor Manfred Lehn for the supervision of her thesis and his suggestions on how to deal with numerous difficulties.
She is grateful to her second advisor Christoph Sorger for proposing her the work on $\theta$--stable $G$--constellations and for the time she spent in Nantes.
In addition, she would like to thank Michel Brion, S\"onke Rollenske, Ziyu Zhang and Markus Zowislok for some fruitful discussions about her work. The second-named author thanks Christian Lehn for useful discussions, and especially for suggesting the correct $\theta$--stability condition in our setting. Both authors gratefully acknowledge the financial support by DAAD and SFB/TR 45.
\end{Ack}

% *** 2 ***************************************************************************************************************************************
\section{$(G,h)$--constellations}\label{chGconst}

To begin we adapt the notion of $G$--constellation, originally introduced by Craw and Ishii in \cite{CI:2004} for finite groups, to the case of infinite reductive groups. In our definition we replace the  isotypic decomposition of the regular representation by an isotypic decomposition given by a prescribed Hilbert function $h$. Further, we adapt Craw and Ishii's notion of $\theta$--stability and $\theta$--semistability and we introduce the moduli functors $\fsMth$ and $\fMth$ of $\theta$--stable and $\theta$--semistable $(G,h)$--con\-stellations, respectively. Then in Section \ref{sfiniteness} we show that $\theta$--semistable $(G,h)$--con\-stellations satisfy a certain finiteness condition. Afterwards, we examine flat families of $(G,h)$--constellations and reduce the verification of the $\theta$--(semi)stability condition to finitely many subsheaves only. The aim is to construct a moduli space of $\theta$--stable $(G,h)$--constellations representing $\fsMth$, which, for a special choice of $\theta$, recovers the invariant Hilbert scheme. Indeed, in Section \ref{HilbgleichMth} we show that if $h(\rho_0) = 1$ and $\theta$ is chosen appropriately, then $\fsMth$ coincides with the invariant Hilbert functor.

\subsection{Definitions}  \label{deff}

We fix an infinite reductive group $G$, an affine $G$--scheme of finite type $X$ and a Hilbert function $h \colon \Irr G \to \natN_0$, where $\Irr G$ denotes the set of isomorphy classes of irreducible representations $\rho\colon G \to Gl(V_{\rho})$. Throughout this article, we consider schemes and algebraic groups over the field of complex numbers $\kxC$. Unless explicitly mentioned, schemes are assumed to be noetherian and $G$--modules are assumed to be rational. 

\begin{definition}\hfill\par
%\vspace{-0.5em}
\begin{enumerate}
\itemsep-0.4ex
 \item Let $R_h := \bigoplus_{\rho \in \Irr G} \kxC^{h(\rho)} \tensor V_{\rho}$ be the $G$--module with multiplicities given by $h$. A $(G,h)$\textit{--constellation} on $X$ is a $G$--equivariant coherent $\ssh_X$--module $\xxF$ such that $H^0(\xxF)$ is isomorphic to $R_h$ as a representation of $G$. In other words $H^0(\xxF)$ is a $H^0(\ssh_X)$--module of finite type equipped with a $G$--module structure such that $H^0(\xxF) \cong R_h$ and
$$\forall g \in G,\; \forall f \in H^0(\ssh_X),\; \forall m \in H^0(\xxF),\ g.(f.m)= (g.f).(g.m).$$ 
 \item Given a scheme $S$, a \textit{family of $(G,h)$--constellations} over $S$ is a coherent sheaf $\xxFf$ on a family of affine $G$--schemes $\mathcal{X}$ over $S$ in the sense of \cite[Definition 1.1]{AB:2005}, i.e.~on a scheme $\mathcal{X}$ equipped with an action of $G$ and an affine $G$--invariant morphism $\mathcal{X} \to S$ of finite type, such that the restrictions $\xxFf(s) = \xxFf|_{\XXX(s)}$ are $(G,h)$--constellations on the fibres $\XXX(s) := \mathcal{X} \times_S \Spec(k(s))$.
\end{enumerate}
\end{definition}

We would like to represent the functor that assigns to a scheme $S$ the set of families of $(G,h)$--constellations on the scheme $X$. In general, the set of $(G,h)$--con\-stellations on $X$ is too large to be parametrized by a scheme. Hence, to construct a moduli space of these objects, we restrict ourselves to $(G,h)$--constellations satisfying a certain stability condition $\theta \in \ratQ^{\Irr G}$ where $\theta_{\rho} < 0$ for only finitely many $\rho \in \Irr G$. This induces a decomposition 
\begin{equation}  \label{def_D}
\Irr G = D_+ \cup D_0 \cup D_- \qquad \text{such that} \qquad
\theta_{\rho}\left\{ \begin{array}{ll} 
> 0, & \rho \in D_+, \\
= 0, & \rho \in D_0, \\
< 0, & \rho \in D_-.
\end{array}\right.
\end{equation}
By the assumption on $\theta$, the set $D_-$ is finite and will always be finite along this paper. To define the $\theta$--stability condition, we first need to associate to $\theta$ a function on the representation monoid $R(G) = \bigoplus_{\rho \in \Irr G}\natN \cdot \rho$ and on the category $\Coh^G(X)$ of $G$--equivariant coherent $\ssh_X$--modules:

\begin{definition}
If $\theta \in \ratQ^{\Irr G}$ with $\theta_{\rho} < 0$ for only finitely many $\rho \in \Irr G$, we define a function $\theta\colon R(G) \to \reR \cup \{\infty\}$ by
$$\theta(W) := \langle\theta,h_{W}\rangle := \sum_{\rho \in \Irr G} \theta_{\rho}\cdot \dim W_{\rho}$$
where $W = \bigoplus_{\rho \in \Irr G}W_{\rho}\tensor V_{\rho}$ is the isotypic decomposition of $W$.

In order to consider $\theta$ as a function $\theta\colon \Coh^G(X) \to \reR \cup \{\infty\}$ we set
$$\theta(\xxF) := \theta(H^0(\xxF)) = \sum_{\rho \in \Irr G} \theta_{\rho}\cdot \dim \xxF_{\rho}$$
where $\xxF_{\rho}:=\xxHom^G(V_\rho,H^0(\xxF))$, and $H^0(\xxF) = \bigoplus_{\rho \in \Irr G}\xxF_{\rho}\tensor V_{\rho}$ is the isotypic decomposition of $H^0(\xxF)$.
In particular, if $\xxF$ is a $(G,h)$--constellation, then we have $\theta(\xxF) = \sum_{\rho \in \Irr G} \theta_{\rho}h(\rho)$.
\end{definition}

A consequence of the finiteness of $D_-$ is that the series $\sum_{\rho \in \Irr G} \theta_{\rho}h(\rho)$ is convergent if and only if the series $\sum_{\rho \in D_+} \theta_{\rho}h(\rho)$ is convergent. 

\begin{definition} 
If a $(G,h)$--constellation $\xxF$ on $X$ is generated by $\bigoplus_{\rho \in D_-} \xxF_{\rho} \tensor V_{\rho}$ as an $\ssh_X$--module, we say that $\xxF$ \textit{is generated in} $D_-$. 
\end{definition} 

We are now in the position to define the stability condition we need on $(G,h)$--constellations:

\begin{definition}
A $(G,h)$--constellation $\xxF$ is called $\theta$\textit{--semistable} if $\theta(\xxF) = 0$ and if for all $G$--equivariant coherent subsheaves $\xxF' \subset \xxF$ generated in $D_-$ we have $\theta(\xxF') \geq 0$. Moreover, $\xxF$ is called $\theta$\textit{--stable} if $\theta(\xxF) = 0$ and if for all non--zero proper $G$--equivariant coherent subsheaves $0 \neq \xxF' \subsetneq \xxF$ generated in $D_-$ we have $\theta(\xxF') > 0$.

For the reader's convenience, we replace the similar conditions for stability and semistability by setting everything concerning semistability in parentheses and we introduce the symbol ``$\sgeq$'': A $(G,h)$--constellation $\xxF$ is called $\theta$\textit{--(semi)stable} if $\theta(\xxF) = 0$ and if for all non--zero proper $G$--equivariant coherent subsheaves $\xxF' \subset \xxF$ generated in $D_-$ we have $\theta(\xxF') \sgeq 0$.
In the same way, ``$\sleq$'' stands for ``$\leq$'' in the case of semistability and ``$<$'' in the case of stability.
\end{definition}

In a former version of this article the definition of $\theta$--(semi)stability required $\theta(\xxF') \sgeq 0$ for all $G$--equivariant coherent subsheaves $\xxF' \subset \xxF$ exactly as in the case of finite groups (\cite[\S 2.1]{CI:2004}). This condition turned out not to work in some proofs so that it had to be modified. The difference lies in those $G$--equivariant coherent sheaves $\xxF$ containing $G$--equivariant coherent subsheaves generated in positive parts. 

\begin{remark} \label{thetaquot}
Every $G$--equivariant subsheaf $\xxF'$ of $\xxF$ induces a $G$--equivariant quotient $\xxF'' := \xxF/\xxF'$ of $\xxF$. Conversely, every $G$--equivariant quotient $\alpha\colon \xxF \twoheadrightarrow \xxF''$ induces a $G$--equivariant subsheaf $\xxF' := \ker \alpha$ of $\xxF$. In both cases the corresponding Hilbert functions satisfy $h_{\xxF'} + h_{\xxF''} = h$, so that $\theta(\xxF) = \theta(\xxF') + \theta(\xxF'')$. Thus a $(G,h)$--constellation $\xxF$ is $\theta$--(semi)stable if and only if $\theta(\xxF) = 0$ and for all non-zero proper $G$--equivariant quotients $\xxF \twoheadrightarrow \xxF''$ whose kernel is generated in $D_-$ we have $\theta(\xxF'') \sleq 0$.
\end{remark}

Since $\theta(\xxF)$ is supposed to be $0$ for any $\theta$---semistable $(G,h)$--constellation $\xxF$, the values of $\theta$ have to be chosen such that $\langle\theta,h\rangle = 0$. In particular, the series $\sum_{\rho \in \Irr G} \theta_{\rho}h(\rho)$ is convergent.

If $\theta = 0$ or at least $\theta_{\rho} = 0$ whenever $h(\rho) \neq 0$, then every $(G,h)$--constellation is $\theta$--semistable, but the only $\theta$--stable $(G,h)$--constellations are those having no \linebreak $G$--equivariant subsheaves generated in $D_-$ different from $0$ and themselves. This case is not of any interest to us. To avoid this, in the following we will always assume that $D_- \cap \supp h$ and $D_+ \cap \supp h$ are non--empty.

Moreover, we will always assume that $\theta_\rho=0$ whenever $h(\rho)=0$. This assumption will simplify the computations in Subsection \ref{comparison} and is not restrictive at all since we always consider the product $\theta_\rho h(\rho)$ or $\theta_\rho h'(\rho)$ (with $h'(\rho) \leq h(\rho)$) in our calculations.

\vspace{1em}

Now we define the moduli functors that we will consider in the following:

\begin{definition}  \label{thetafunctors}
The moduli functor of $\theta$--semistable $(G,h)$--constellations on $X$ is 
\begin{align*}
&\fMth\colon \text{(Sch/$\kxC$)}^{\text{op}} \to \text{(Set)} \\
&\;S \mapsto \{\xxFf \text{ an } S\text{--flat family of $\theta$--semistable $(G,h)$--constellations on } X \times S\}/_{\cong}, \\
&\;(f\colon S' \to S) \mapsto \big(\fMth(S) \to \fMth(S'), \xxFf \mapsto (\xid_X\times f)^*\xxFf\big).\\
\intertext{The moduli functor of $\theta$--stable $(G,h)$--constellations on $X$ is}
&\fsMth\colon \text{(Sch/$\kxC$)}^{\text{op}} \to \text{(Set)} \\
&\;S \mapsto \{\xxFf \text{ an } S\text{--flat family of $\theta$--stable $(G,h)$--constellations on } X \times S\}/_{\cong}, \\
&\;(f\colon S' \to S) \mapsto \big(\fsMth(S) \to \fsMth(S'), \xxFf \mapsto (\xid_X\times f)^*\xxFf\big).
\end{align*}
\end{definition}

\subsection{Finiteness} \label{sfiniteness}

Our strategy to construct the moduli space $\sMth$ of $\theta$--stable $(G,h)$--constel-lations is to show that all $\theta$--stable $(G,h)$--constellations are quotients of a certain $G$--equivariant coherent $\ssh_X$--module $\xxH$ and to obtain our moduli space by considering the invariant Quot scheme $\xxGQuotH$ and its GIT--quotient.

In Theorem \ref{finiteness} we will see that $\theta$--stable $(G,h)$--constellations are controlled by their isotypic components in $D_-=\{\rho \in \Irr G \ | \ \theta_\rho<0\}$. As $D_-$ is finite this means that all $\theta$--stable $(G,h)$--constellations are generated by finitely many irreducible representations.

Let $\xxF$ be a $\theta$--(semi)stable $(G,h)$--constellation and $\xxF'$ a $G$--equivariant coherent subsheaf of $\xxF$ generated in $D_-$. Let $H^0(\xxF') \cong \bigoplus_{\rho \in \Irr G} \xxF'_{\rho} \tensor V_{\rho}$ be the isotypic decomposition of its global sections. Then we have $h'(\rho) := \dim \xxF'_{\rho} \leq h(\rho)$ for every $\rho \in \Irr G$. Since $D_-$ is finite, $\theta(\xxF')$ is also a convergent series and we have
$$\theta(\xxF') = \sum_{\rho \in \Irr G} \theta_{\rho}h'(\rho)
= \underbrace{\sum_{\rho \in D_-} \underbrace{\theta_{\rho}}_{<0}\underbrace{h'(\rho)}_{\geq 0}}_{\leq 0} + \underbrace{\sum_{\rho \in D_+} \underbrace{\theta_{\rho}}_{>0}\underbrace{h'(\rho)}_{\geq 0}}_{\geq 0} \quad \overset{!}{\sgeq} \quad 0.$$

As a philosophy, if $\xxF$ is to be $\theta$--(semi)stable, the values $h'(\rho)$ should be as large as possible in $D_+$ and as small as possible in $D_-$. This means that all subsheaves of $\xxF$ generated in $D_-$ should be similar to $\xxF$ in positive parts and they should nearly vanish in negative parts. 
%We can thus expect that the subsheaf of $\xxF$ generated by its summands in $D_-$ will almost coincide with $\xxF$.  
%In other words, the most destabilizing subsheaf of $\xxF$ is the subsheaf of $\xxF$ generated by its summands in $D_-$.

We have the following finiteness result:

\begin{thm}\label{finiteness}
Let $\theta \in \ratQ^{\Irr G}$ such that $\theta_{\rho} < 0$ for only finitely many $\rho \in \Irr G$. If $\xxF$ is a $\theta$--stable $(G,h)$--constellation on $X$, then $\xxF$ is generated in $D_-$. Moreover, if $D_0 = \{ \rho \in \Irr G \ |\ h(\rho)=0\}$, the same holds for any $\theta$--semistable $(G,h)$--constel\-lation on $X$.
\end{thm}

\begin{proof}
Consider the $\ssh_X$--submodule $\xxF'$ of $\xxF$ generated by $\bigoplus_{\rho \in D_-} \xxF_{\rho} \tensor V_{\rho}$. Then we have:
$$\begin{array}{rll}
h'(\rho) = h(\rho) &\text{ for } \;\rho \in D_-, \\
h'(\rho) \leq h(\rho) &\text{ for } \;\rho \in D_+ \cup D_0.
\end{array}$$
This implies
$$\theta(\xxF') = \sum_{\rho \in D_-} \theta_{\rho}h'(\rho) +  \sum_{\rho \in D_+} \theta_{\rho}h'(\rho)
\leq \sum_{\rho \in D_-} \theta_{\rho}h(\rho) +  \sum_{\rho \in D_+} \theta_{\rho}h(\rho) = \theta(\xxF) = 0.$$
If $\xxF$ is $\theta$--stable this means that $\xxF' = \xxF$, because otherwise $\xxF'$ would destabilize $\xxF$. If $\xxF$ is $\theta$--semistable we obtain $\theta(\xxF')=0$ and thus $h'(\rho)=h(\rho)$ for every $\rho \in D_+$. As $D_0 = \{ \rho \in \Irr G \ |\ h(\rho)=0\}$ we obtain that $h'=h$, that is $\xxF'=\xxF$.
This shows that every $\theta$--(semi)stable sheaf $\xxF$ is generated by $\bigoplus_{\rho \in D_-} \xxF_{\rho} \tensor V_{\rho}$.
\end{proof}

This finiteness result causes us to define the following $G$--equivariant free \linebreak $\ssh_X$--module of finite rank:
\begin{equation}\label{defH}\xxH := \left( \bigoplus_{\rho \in D_-} \kxC^{h(\rho)} \tensor V_{\rho} \right) \tensor \ssh_X \cong \ssh_X^{\sum_{\rho \in D_-} h(\rho)\dim V_{\rho}}.\end{equation}
Then by Theorem \ref{finiteness} it follows that every $\theta$--(semi)stable $(G,h)$--constellation can be obtained as a quotient of $\xxH$ if $D_-$ is finite (and if $D_0=\{ \rho \in \Irr G \ |\ h(\rho)=0\}$). We will establish this in more detail in Section \ref{quotconst}. Consequently, we may consider $\xxGQuotH$ to construct the moduli space of $\theta$--(semi)stable $(G,h)$--constellations.

Another consequence of the consideration of $D_-$ is that $\theta$--(semi)stability can be proven by checking finitely many subsheaves only, as the following sequence of results shows.

%\begin{prop}\label{Dminusgenuegt}
%Let $\xxF$ be a $(G,h)$--constellation with $\theta(\xxF) = 0$. If for every non-zero proper $G$--equivariant subsheaf $\widetilde \xxF \subset \xxF$ generated in $D_-$ we have $\theta(\widetilde \xxF) \sgeq 0$, then $\xxF$ is $\theta$--(semi)stable.
%\end{prop}
%\begin{proof}
%Assume that $\theta(\widetilde \xxF) \sgeq 0$ for every proper $G$--equivariant subsheaf $\widetilde \xxF \subset \xxF$ generated in $D_-$ and let $\xxF'$ be a $G$--equivariant subsheaf of $\xxF$. Consider the subsheaf $\xxF^*$ of $\xxF'$ generated by the $\xxF'_{\rho}$, $\rho \in D_-$, so that $h^*(\rho) := \dim \xxF_{\rho}^* = h'(\rho)$ for $\rho \in D_-$ and $h^*(\rho) \leq h'(\rho)$ for $\rho \in \Irr G \setminus D_-$. Since $\xxF^*$ is generated in $D_-$, we have
%\begin{align*}
%\theta(\xxF') &= \sum_{\rho \in D_-}\theta_{\rho}h'(\rho) \; + \sum_{\rho \in \Irr G \setminus D_-}\underbrace{\theta_{\rho}}_{\geq 0}h'(\rho) \\
%&\geq \sum_{\rho \in D_-}\theta_{\rho}h^*(\rho) \; + \sum_{\rho \in \Irr G \setminus D_-}\theta_{\rho}h^*(\rho) = \theta(\xxF^*) \sgeq 0.
%\end{align*} 
%\end{proof}

\begin{lem}\label{bounded}
The family of pairs 
\begin{equation}\label{pairs}
\left\{(\xxF,\xxF') \, \middle| \, \begin{array}{ll}\xxF \text{ a $(G,h)$--constellation generated in $D_-$,} \\ \xxF' \subset \xxF \text{ a $G$--equivariant coherent subsheaf generated in $D_-$}\end{array}\right\}
\end{equation}
is bounded, i.e.~there is a scheme $Z$ of finite type, a $G$--equivariant coherent sheaf of $\ssh_{X\times Z}$--modules $\xxFf$ and a $G$--equivariant coherent subsheaf $\xxFf'$ of $\xxFf$ such that the family \eqref{pairs} is a subset of $\{(\xxFf|_{X\times \Spec(k(z))}, \xxFf'|_{X\times \Spec(k(z))}) \mid z \text{ a closed point in } Z\}$.
\end{lem}

\begin{proof}
The set of $(G,h)$--constellations $\xxF$ generated in $D_-$ is parametrized by the quasiprojective scheme $\xxGQuotH$. For a fixed $\xxF$ the subsheaves $\xxF' \subset \xxF$ generated in $D_-$ are determined by the choice of subspaces $\xxF'_{\rho} \subset \xxF_{\rho}$ for $\rho \in D_-$. %Hence they are parametrized by finite subsets of $\prod_{\rho \in D_-} \coprod_{k = 0}^{h(\rho)} \xxGrass{k}{\kxC^{h(\rho)}}$. 
Hence the set \eqref{pairs} is parametrized by a subset of the scheme 
$$Z := \xxGQuotH \times \prod_{\rho \in D_-} \coprod_{k = 0}^{h(\rho)} \xxGrass{k}{\kxC^{h(\rho)}}.$$ 
The scheme $Z$ is quasiprojective, hence of finite type, and the family \eqref{pairs} is bounded by the universal family of its functor of points.
\end{proof}

\begin{prop}\label{endlvieleHilbertfkt}   
There is a finite set of Hilbert functions $\{ h_1,\ldots, h_n\}$ such that for any $(G,h)$--constellation $\xxF$ and any $G$--equivariant coherent subsheaf $\xxF' \subset \xxF$, both generated in $D_-$, the Hilbert function $h'$ of $\xxF'$ is one of the $h_1,\ldots,h_n$.
\end{prop}

\begin{proof}
Lemma \ref{bounded} says that the family of pairs $(\xxF,\xxF')$ with $\xxF$ a $(G,h)$--constella\-tion and $\xxF'$ a $G$--equivariant coherent subsheaf of $\xxF$, both generated in $D_-$, is bounded by a pair of coherent sheaves $(\xxFf,\xxFf')$ on $X \times Z$, where $Z$ is a scheme of finite type. 
The family $(\xxFf,\xxFf')$ is not necessarily flat on $Z$, but we can use \cite[Lemme 3.4]{FGAIV:1961} to obtain a flattening stratification of $Z$, that is a finite decomposition $Z = \coprod_{i = 1}^n Z_i$ of $Z$ into a disjoint union of connected and locally closed subschemes $Z_i \subset Z$ such that $(\xxFf|_{Z_i},\xxFf'|_{Z_i})$ is a flat family on $Z_i$. Then for all $z \in Z_i$ the fibres $\xxFf'(z)$ have the same Hilbert function $h_i$.
%Lemme 3.4 in \cite{FGAIV:1961} is only formulated in the case where $\ssh_{X \times Z}$ and $(\xxFf,\xxFf')$ are graded over $\natN_0$, $\ssh_{X \times Z}$ is generated by $(\ssh_{X \times Z})_1$ and $h$ is a polynomial. We reduce our situation to this setting as follows:
%Define a map $a\colon \Irr G \cong \natN_0^{\rk G} \to \natN_0$ via $\rho = \sum_{\rho_i \in \Irr G}n_i\rho_i \mapsto \sum_{\rho \in \Irr G} n_i$, where all but finitely many $n_i$ vanish. Then $\ssh_{X \times Z}$ is graded over $\natN_0$ with $(\ssh_{X \times Z})_n = \bigoplus_{a(\rho) = n}(\ssh_{X \times Z})_{\rho}$. The same holds for $\xxFf$ and $\xxFf'$. The function $p\colon \natN_0 \to \natN_0$, $p(n) = \sum_{a(\rho) = n}h(\rho)$ describes the rank of the $\xxFf_n$ and analogously we have $p'$ for $\xxFf'$.
%Further, there is a degree $d$ such that the ring $\ssh_{X \times Z}^{(d)} := \bigoplus_{n \in \natN_0}(\ssh_{X \times Z})_{nd}$ is generated by $(\ssh_{X \times Z}^{(d)})_1 = (\ssh_{X \times Z})_{d}$. For $i = 1, \hdots, d-1$ set $\xxFf^i := \bigoplus_{n \in \natN_0}\xxFf_{i + nd}$, so that $\xxFf = \xxFf^0 \oplus \hdots \oplus \xxFf^{d-1}$. Then all the $\xxFf^i$ are $\ssh_{X \times Z}^{(d)}$--modules and each corresponding function $p^i$ with $p^i(n) = \rk \xxFf^i_n$ is a polynomial. By \cite[Lemme 3.4]{FGAIV:1961} we find a flattening stratification for each $\xxFf^i$. In the same way we obtain a flattening stratification for the $(\xxFf')^i$. Their common refinement yields a flattening stratification for $(\xxFf, \xxFf')$. 
\end{proof}

\begin{cor}
Let $\theta \in \ratQ^{\Irr G}$ such that $\theta_{\rho} < 0$ for only finitely many $\rho \in \Irr G$, and let $\xxF$ be a $(G,h)$--constellation generated in $D_-$ with $\theta(\xxF) = 0$. With the notation of Proposition \ref{endlvieleHilbertfkt} suppose that for every $i = 1,\hdots,n$ with $h_i$ actually occurring as a Hilbert function of some non-zero proper $G$--equivariant subsheaf of $\xxF$ generated in $D_-$, we have $\langle\theta,h_i\rangle \sgeq 0$. Then $\xxF$ is $\theta$--(semi)stable.
\end{cor}

\subsection{The invariant Hilbert scheme as a moduli space of $(G,h)$--constel\-lations} \label{HilbgleichMth}

Let us suppose that the Hilbert function $h$ satisfies $h(\rho_0)=1$, where $\rho_0$ is the trivial representation. For recovering the invariant Hilbert functor (cf.~\cite[Definition 1.5]{AB:2005} or \cite[Definition 2.1]{Bec:2011}) and the invariant Hilbert scheme, $\theta$ must satisfy an extra condition: 
\begin{prop}
If $h(\rho_0)=1$ and $\theta$ is chosen such that $D_- = \{\rho_0\}$, then the moduli functor of $\theta$--stable $(G,h)$--constellations coincides with the invariant Hilbert functor: $$\fsMth = \finvHilb{h}(X).$$
\end{prop}

\begin{proof} Let $S$ be a noetherian scheme over $\kxC$, $s \in S$ a point and $\xxF = \xxFf(s)$ a fibre of a flat family $\xxFf$ of $\theta$-stable $(G,h)$--constellations on $X \times S$. Theorem \ref{finiteness} and the condition $D_- = \{\rho_0\}$ imply that the $\ssh_X$--module generated by $V_{\rho_0}$ is $\xxF$, i.e.~$\xxF$ is cyclic and therefore it is isomorphic to a quotient of $\ssh_X$. This means $\xxF \cong \ssh_{Z_s}$ for some ${Z_s} \in \invHilb{G}{h}{X}$ and setting $\xxZ = \{(Z_s,s)\mid s \in S\} \in \finvHilb{h}(X)(S)$ we obtain $\xxFf \cong \ssh_{\xxZ}$.

%means $\theta_{\rho_0}h(\rho_0) = -\sum\limits_{\rho \in D_+} \underbrace{\theta_{\rho}}_{>0}\underbrace{h(\rho)}_{>0} < 0$. For any non-zero proper $G$--equivariant subsheaf $\xxF' \subset \xxF$ we have $\theta(\xxF') = \sum\limits_{\rho \in \Irr G} \theta_{\rho}\underbrace{h'(\rho)}_{\leq h(\rho)}$. Taking into account that $h(\rho_0) = 1$ there are two cases for $h'(\rho_0)$:
%\begin{itemize}\renewcommand{\labelitemi}{\textbullet}
% \item $h'(\rho_0) = 1 = h(\rho_0)$: In this case
%$$\theta(\xxF') = \theta_{\rho_0}\cdot 1 \; + \sum_{\begin{smallmatrix}\rho \in \Irr G\\ \rho \neq \rho_0 \end{smallmatrix}} \theta_{\rho}h'(\rho) 
%= \sum_{\begin{smallmatrix}\rho \in \Irr G\\ \rho \neq \rho_0 \end{smallmatrix}} \underbrace{\theta_{\rho}}_{\geq0}(\underbrace{h'(\rho) - h(\rho)}_{\leq 0}) \leq 0,$$
%so for stable $\xxF$ this case cannot occur.
%\item Hence for stable $\xxF$ we have $h'(\rho_0) = 0$, so that no proper subsheaf of $\xxF$ contains $V_{\rho_0}$. Thus the $\ssh_X$--module generated by $V_{\rho_0}$ is $\xxF$, i.e.~$\xxF$ is cyclic and therefore it is isomorphic to a quotient of $\ssh_X$. This means $\xxF \cong \ssh_{Z_s}$ for some ${Z_s} \in \invHilb{G}{h}{X}$ and setting $\xxZ = \{(Z_s,s)\mid s \in S\} \in \finvHilb{h}(X)(S)$ we obtain $\xxFf \cong \ssh_{\xxZ}$.
%\end{itemize}
Conversely, consider an element $\xxZ \in \finvHilb{h}(X)(S)$. Every fibre $\ssh_{\xxZ}(s)$ of its structure sheaf is generated by the image of $1 \in \ssh_X$, which is an invariant. Therefore, every non-zero proper $G$--equivariant subsheaf $\xxF'$ of $\ssh_{\xxZ}(s)$ satisfies $h'(\rho_0) = 0$ and hence is generated in $D_0 \cup D_+$. So $\ssh_{\xxZ}(s)$ is $\theta$--stable for every $s \in S$, which means $\ssh_{\xxZ} \in \fsMth(S)$.
%$\theta(\xxF') > 0$ (since $D_0=\emptyset$). So $\ssh_{\xxZ}(s)$ is $\theta$--stable for every $s \in S$, which means $\ssh_{\xxZ} \in \fsMth(S)$.
\end{proof}

\begin{cor}\label{recoverGHilb}
If $h(\rho_0)=1$ and $\theta$ is chosen such that $D_- = \{\rho_0\}$, then the functor $\fsMth$ is representable and the moduli space of $\theta$--stable $(G,h)$--constellations is $\sMth = \invHilb{G}{h}{X}$.
\end{cor}

% *** 3 ***************************************************************************************************************************************
\section[GIT of the invariant Quot scheme]{Geometric Invariant Theory of the invariant Quot scheme}\label{chQuotGIT}

In the last section we have shown that every $\theta$--(semi)stable $(G,h)$--constellation is a quotient of $\xxH := \bigoplus_{\rho \in D_-} \kxC^{h(\rho)} \tensor V_{\rho} \tensor \ssh_X$. Now we consider the invariant Quot scheme $\xxGQuotH$ parametrizing all $G$--equivariant quotient maps $[q\colon \xxH \twoheadrightarrow \xxF]$, where $\xxF$ is a $G$--equivariant coherent $\ssh_X$--module whose module of global sections is isomorphic to $R_h := \bigoplus_{\rho \in \Irr G} V_{\rho}^{\oplus h(\rho)}$. 
In Subsection \ref{invQuot} we consider an embedding of the invariant Quot scheme into a product of Grassmannians. This equips $\xxGQuotH$ with an ample line bundle $\xLl$.
Thereafter we discuss the geometric invariant theory (GIT) of $\xxGQuotH$ in order to obtain a categorical quotient $\xxGQuotH^{ss} \red _{\xLl_{\chi}} \Gamma$ of GIT--semistable quotients and a geometric quotient $\xxGQuotH^{s} \red _{\xLl_{\chi}} \Gamma = \xxGQuotH^{s}/\Gamma$ of GIT--stable quotients. 
The \linebreak $\theta$--stable $(G,h)$--constellations identify with elements of the latter and it will be our candidate for the moduli space of $\theta$--stable $(G,h)$--constellations. Here, $\Gamma$ denotes the gauge group of $\xxH$ and $\xLl_{\chi}$ is the ample line bundle $\xLl$ with linearization depending on the choice of a character $\chi$ of $\Gamma$. We describe these parameters in Subsection \ref{GITparameters}. Afterwards, in Subsection \ref{onePS} we examine $1$--parameter subgroups of $\Gamma$ and establish their description via filtrations of the vector space $\bigoplus_{\rho \in D_-} \kxC^{h(\rho)}$ in order to apply Mumford's numerical criterion for GIT--(semi)stability in Subsection \ref{GITstab1}. Out of this we eventually establish a condition for GIT--(semi)stability by considering subspaces of $\bigoplus_{\rho \in D_-} \kxC^{h(\rho)}$ instead of filtrations. This condition will be used to compare GIT--(semi)stability to $\theta$--(semi)stability in Section \ref{chcorres}.

\subsection{Embedding of the invariant Quot scheme}\label{invQuot}   

Let $\xxH$ be any $G$--equivariant coherent $\ssh_X$--module with isotypic decomposition $H^0(\xxH) = \bigoplus_{\rho \in \Irr G} \xxH_{\rho} \tensor V_{\rho}$ and $h\colon \Irr G \to \natN_0$ a Hilbert function. Then we consider the invariant Quot scheme $\xxGQuotH$ as constructed in \cite{Jan:2006}.
Before we address ourselves to the geometric invariant theory of the invariant Quot scheme, we consider the embedding of $\xxGQuotH$ into a finite product $\prod_{\sigma \in D}\qGrass{H_{\sigma}}{h(\sigma)}$ of Grassmannians generalizing the embedding of the invariant Hilbert scheme \cite[Section 4.2]{Bec:2011}.

The next result follows directly from the construction of the invariant Quot scheme by Jansou \cite[\S1.2 and \S1.3]{Jan:2006}.

\begin{prop}  \label{embedding_etha}
There exists a finite subset $D \subset \Irr G$ and, for each $\rho \in D$, a finite dimensional vector space $H_\rho$ together with a surjection of $\ssh_{X\red G}$--modules $\ssh_{X\red G} \otimes_\kxC H_\rho \twoheadrightarrow \xxH_\rho$ such that there is a locally closed immersion
\begin{equation}\label{Quotimmersion}
\eta\colon \xxGQuotH \longhookrightarrow \prod\limits_{\rho \in D} \qGrass{H_{\rho}}{h(\rho)}.
\end{equation}
\end{prop}

\begin{proof}
Let us explain how the construction of the invariant Quot scheme by Jansou gives the desired embedding into a product of Grassmannians. 

Fix a Borel subgroup $B \subset G$, and denote by $T \subset B$ a maximal torus and by $U$ the unipotent radical of $B$. By \cite[Proposition 1.10]{Jan:2006}, there is a closed immersion $\xxGQuotH \hookrightarrow \operatorname{Quot}^T(\xxH^U,h')$, where $h': \Irr T \to \natN_0$ coincides with $h$ on $\Irr G$ and is $0$ elsewhere, and $\xxH^U$ is the $T$--equivariant coherent $\ssh_{X\red U}$--module obtained by taking the $U$--invariants of $\xxH$. 

In this setting let $E$ be a finite dimensional $T$--module such that $X\red U$ identifies with a $T$--stable closed subscheme of $E$, let $e_1,\ldots,e_r$ be a system of generators of the $\ssh_{X\red U}$--module $\xxH^U$ formed by weight vectors for the action of $T$, and define \linebreak $\xMM:= \bigoplus_{i=1}^{r} \ssh_E e_i$, which is a $T$--equivariant coherent $\ssh_E$--module. By \cite[Lemme 1.9]{Jan:2006}, there exists a closed immersion $\operatorname{Quot}^T(\xxH^U,h') \hookrightarrow \operatorname{Quot}^T(\xMM,h')$.  Besides, by definition of $\xMM$, there is a surjective morphism of $T$--equivariant $\ssh_E$--modules $\xMM \twoheadrightarrow \xxH^U$, where the structure of $\ssh_E$--module for $\xxH^U$ comes from the surjection $\ssh_E \twoheadrightarrow \ssh_{X\red U}$. In particular, for each $\rho \in \Irr T$, we have a surjection of $\ssh_{E\red T}$--modules $\xMM_\rho \twoheadrightarrow \xxH_{\rho}^{U}=\xxH_\rho$. 

Now by \cite[Proposition 1.6, Lemme 1.7, and Proposition 1.8]{Jan:2006}, there exists a finite subset $D \subset \Irr T$ and, for each $\rho \in D$, a finite dimensional subspace $H_\rho \subset \xMM_\rho$ such that $\operatorname{Quot}^T(\xMM,h')$ identifies with a locally closed subscheme of the product $\prod_{\rho \in D} \qGrass{H_{\rho}}{h'(\rho)}$ and such that $H_\rho$ generates $\xMM_\rho$ as an $\ssh_{E\red T}$--module. We thus obtain a surjection of $\ssh_{E\red T}$--modules $\ssh_{E\red T} \otimes_\kxC H_\rho \twoheadrightarrow \xxH_\rho$, but any element in the kernel of $\ssh_E \twoheadrightarrow \ssh_{X\red U}$ acts trivially on $\xxH^U$, whence a surjection of $\ssh_{X\red G}$--modules $\ssh_{X\red G} \otimes_\kxC H_\rho \twoheadrightarrow \xxH_\rho$. 

Finally, as $h'(\rho)=0$ for every $\rho \in \Irr T \backslash \Irr G$, we can assume that $D \subset \Irr G$ and the theorem is proven. 
\end{proof}

\begin{remark}  \label{extension D}
The existence of the set $D$ in Proposition \ref{embedding_etha} is given by \cite[Lemme 1.7]{Jan:2006} and one easily checks that any finite subset $D'$ of $\Irr G$ containing $D$ also provides an embedding of the invariant Quot scheme.  
\end{remark}

\subsection{The parameters needed for GIT}\label{GITparameters}

We fix $\theta \in \ratQ^{\Irr G}$ such that the set $D_-$ defined by \eqref{def_D} is finite, and let $\xxH$ be as defined in \eqref{defH}.
In this subsection we introduce a group action on the invariant Quot scheme of $\xxH$, for which we want to obtain the GIT--quotient. In order to determine this quotient, we need to find an ample line bundle on $\xxGQuotH$, which can be linearized with respect to the group action. The linearization depends on a character of the group.

In the definition of $\xxH$, we write $A_{\rho} := \kxC^{h(\rho)}$, i.e.~$\xxH := \bigoplus_{\rho \in D_-} A_{\rho} \tensor V_{\rho} \tensor \ssh_X$.
For every $[q \colon \xxH \twoheadrightarrow \xxF] \in \xxGQuotH$, the sheaf $\xxF = q(\xxH)$ is thus generated by the finitely many components $q(A_{\rho} \tensor V_{\rho} \tensor 1)$, $\rho \in D_-$, as an $\ssh_X$--module. 

\subsubsection{The line bundle $\xLl$ and the weights $\kappa$}  \label{linebundleL}

In the last subsection we have seen that there is a finite subset $D \subset \Irr G$ and an embedding $\eta$ of $\xxGQuotH$ into a product of Grassmannians $\prod_{\sigma \in D} \qGrass{H_{\sigma}}{h(\sigma)}$, where the $H_{\sigma}$ are the finite dimensional vector spaces given by Proposition \ref{embedding_etha}. 
Composing $\eta$ with the Pl\"ucker embedding $\pi_{\sigma}$ for every occurring Grassmannian we have
\begin{equation}\label{Quotembedding}\xxGQuotH \overset{\eta}{\longhookrightarrow} \prod_{\sigma \in D} \qGrass{H_{\sigma}}{h(\sigma)} \overset{(\pi_{\sigma})_{\sigma}}{\longhookrightarrow} \prod_{\sigma \in D} \prP(\Lambda^{h(\sigma)}H_{\sigma}).\end{equation}

As mentioned in Remark \ref{extension D}, for any set containing $D$ we again obtain an embedding. Adding further representations if necessary, we will always assume $D_- \subset D$. 

In the following discussion of the geometric invariant theory, different choices of $D$ lead to different notions of GIT--(semi)stability. We will take advantage of the variation of $D$ and the corresponding stability condition in Section \ref{comparison}.

For every choice of $\kappa \in \natN_{0}^D$, the ample line bundles $\ssh_{\sigma}(1)$ on $\prP(\Lambda^{h(\sigma)}H_{\sigma})$ give a line bundle $\bigotimes_{\sigma \in D} (\pi_{\sigma}^*\ssh_{\sigma}(1))^{\kappa_{\sigma}} = \bigotimes_{\sigma \in D} (\det \mathcal{W}_{\sigma})^{\kappa_{\sigma}}$ on the product of the Grassmannians, where $\mathcal{W}_{\sigma}$ denotes the universal family of $\qGrass{H_{\sigma}}{h(\sigma)}$. It is ample if $\kappa_{\sigma} \geq 1$ for every $\sigma \in D$. This in turn induces an ample line bundle
\begin{equation}\label{linebdl}\xLl = \eta^*\bigotimes_{\sigma \in D} (\pi_{\sigma}^*\ssh_{\sigma}(1))^{\kappa_{\sigma}} = \bigotimes_{\sigma \in D} (\det \mathcal{U}_{\sigma})^{\kappa_{\sigma}}\end{equation}
on $\xxGQuotH$, where $p_*\mathcal{U} = \bigoplus_{\sigma \in \Irr G} \mathcal{U}_{\sigma} \tensor V_{\sigma}$ is the isotypic decomposition of the universal quotient $[\pi^*\xxH \twoheadrightarrow \mathcal{U}]$ on $X \times \xxGQuotH$. Here, we denote by $\pi\colon X \times \xxGQuotH \to X$ and $p\colon X \times \xxGQuotH \to \xxGQuotH$ the projections.

\begin{remark}
In Section \ref{comparison} we will also consider $\xLl$ with weights $\kappa_{\sigma} \in \ratQ_{>0}$. To give this a meaning, let $k$ be the common denominator of all the $\kappa_{\sigma}$, $\sigma \in D$. Then we have $k\kappa_{\sigma} \in \natN$ for all $\sigma \in D$ and $\xLl^k$ is an ample line bundle on $\xxGQuotH$, which defines an embedding as above.
\end{remark}

\subsubsection{The gauge group $\Gamma$ and the character $\chi$}\label{parameters}

In order to give concrete surjections $\xxH \twoheadrightarrow \xxF$ rather than only coherent $\ssh_X$--modules $\xxF$ which are quotients of $\xxH$, we have to choose a map $A_{\rho} \to \xxF_{\rho}$ for every $\rho \in D_-$. In order to obtain a moduli space parametrizing sheaves $\xxF$ independent of this choice, we need to consider the natural action of the gauge group 
$\Gamma':= \prod_{\rho \in D_-} Gl(A_{\rho})$ on $\xxH$ by multiplication from the left on the constituent components. 

This action induces a natural action on $\xxGQuotH$ from the right: \\
Let $\gamma = (\gamma_{\rho})_{\rho \in D_-} \in \Gamma'$ and $[q\colon \xxH \twoheadrightarrow \xxF] \in \xxGQuotH$. Then $[q]\cdot\gamma$ is the map
$$[q] \cdot \gamma\colon \xxH \twoheadrightarrow \xxF, \ \; a_{\rho} \otimes v_{\rho} \otimes f \mapsto q(\gamma_{\rho}a_{\rho} \otimes v_{\rho} \otimes f).$$
Since the subgroup of scalar matrices $K:=\{ \prod_{\rho \in D_-} \alpha  \mathrm{Id}_{A_\rho}; \alpha \in \kxC^* \} \cong \kxC^*$ acts trivially on the invariant Quot scheme, we actually consider the action (with finite stabilizers) of the subgroup
\begin{equation} \label{def Gamma}
\Gamma:= \left \{ (\gamma_\rho)_{\rho \in D_-} \in \prod_{\rho \in D_-} Gl(A_\rho) \; \middle| \; \prod_{\rho \in D_-} \det(\gamma_\rho)=1 \right \}.
\end{equation}
Further, the action of $\Gamma$ induces a natural linearization on some power $\xLl^k$ of $\xLl$ (compare to the remark after Lemma 4.3.2 in \cite{HL:2010}). Replacing $\kappa_{\sigma}$ by $k\kappa_{\sigma}$ for every $\sigma \in D$, we can assume that $\xLl$ itself carries a $\Gamma$--linearization. 
Additionally, we can twist this linearization with respect to a character $\chi$ of $\Gamma$, where $\chi(\gamma) = \prod_{\rho \in D_-} \det(\gamma_\rho)^{\chi_\rho}$ and $(\chi_\rho)_{\rho \in D_-} \in \intZ^{D_-}$ are chosen such that $\sum_{\rho \in D_-}\chi_{\rho}h(\rho) = 0$, i.e. $\chi$ restricts to the trivial character on $K \cap \Gamma$. This last condition will be useful to obtain equivalence \eqref{onestepstab}. We write $\xLl_{\chi}$ for the line bundle $\xLl$ equipped with the linearization twisted by the character $\chi$.

\subsection{One--parameter subgroups and filtrations}\label{onePS}

For the construction of the GIT--quo\-tient, we examine $1$--parameter subgroups of $\Gamma$ in order to apply Mumford's numerical criterion and hence deduce a condition for GIT--(semi)stability. 
Let $[q\colon \xxH \twoheadrightarrow \xxF] \in \xxGQuotH$ and $\lambda\colon \kxC^* \to \Gamma$ be a $1$--parameter subgroup. 

Then $\lambda$ induces a grading and a descending filtration on $A := \bigoplus_{\rho \in D_-}A_{\rho}$, so that for every $\rho \in D_-$ we have
$$A_{\rho} = \bigoplus_{n \in \intZ} A_{\rho}^n, \qquad A_{\rho}^{\geq n} = \bigoplus_{m \geq n} A_\rho^m,$$
where $A_{\rho}^n = \{a \in A_{\rho} \mid \lambda(t)\cdot a = t^n a\}$ is the subspace of $A_{\rho}$ on which $\lambda$ acts with weight $n$. 

Let us note that $\lambda$ induces a $1$--parameter subgroup of $\Gamma'$ obtained by composing $\lambda$ with the inclusion $\Gamma \subset \Gamma'$. In fact, denoting $A^m:=\bigoplus_{\rho \in D_-} A_{\rho}^{m}$, one easily checks that a $1$--parameter subgroup of $\Gamma'$ factors through a $1$--parameter subgroup of $\Gamma$ if and only if 
\begin{equation}  \label{sumAm}
\sum_{m \in \intZ} m \cdot \dim(A^m)=0.
\end{equation}  

The grading of $A$ induces a grading 
$$ \xxH = \bigoplus_{n \in \intZ} \xxH^n, \quad \text{ where}\quad \xxH^n = \bigoplus_{\rho \in D_-} A_{\rho}^n \tensor V_{\rho} \tensor \ssh_X,$$
and the corresponding filtration is
$$\xxH^{\geq n} = \bigoplus_{m \geq n} \xxH^m  = \bigoplus_{\rho \in D_-} A_{\rho}^{\geq n} \tensor V_{\rho} \tensor \ssh_X.$$ 
This in turn induces a filtration of $\xxF$ by
$$\xxF^{\geq n} := q(\xxH^{\geq n}),$$
and we define graded pieces
$$\xxF^{[n]} := \xxF^{\geq n}/\xxF^{\geq n+1}.$$

\begin{remark}\label{filtrfin}
As $A_\rho$ is a finite dimensional vector space, only finitely many $A_\rho^n$ are non-zero for every $\rho \in D_-$, so the same holds for $\xxH^n$ and $\xxF^{[n]}$. Further, only finitely many $\xxH^{\geq n}$ and $\xxF^{\geq n}$ are different from $0$ or $\xxH$, respectively from $0$ or $\xxF$.
\end{remark}

The graded object corresponding to the filtration of $\xxF$ is 
$$\overline{\xxF} :=  \bigoplus_{n \in \intZ} \xxF^{[n]} = \bigoplus_{n \in \intZ} \xxF^{\geq n}/\xxF^{\geq n+1}.$$
For the sheaves of covariants of $\overline{\xxF}$ we have $\overline{\xxF}_{\sigma} = \bigoplus_{n \in \intZ} \xxF^{[n]}_{\sigma}$ for every $\sigma \in \Irr G$.
Since $G$ is reductive, the sequences
$$0 \to \xxF_{\sigma}^{\geq n + 1} \to \xxF_{\sigma}^{\geq n} \to \xxF_{\sigma}^{[n]} \to 0$$
are exact for every $\sigma \in \Irr G$, $n \in \intZ$, so that $\dim \xxF_{\sigma}^{[n]} = \dim \xxF_{\sigma}^{\geq n}  - \dim \xxF_{\sigma}^{\geq n + 1}$. 
Let $M,N \in \intZ$ such that $\dim \xxF_{\sigma}^{[n]} = 0$ for every $n > M$, $n < -N$. Then $\xxF_{\sigma}^{\geq -N} = \xxF_{\sigma}$ and $\xxF_{\sigma}^{\geq M + 1} = 0$ and we have
\begin{align*}
\dim \overline{\xxF}_{\sigma} &= \sum_{n \in \intZ} \dim \xxF_{\sigma}^{[n]} = \sum_{n = -N}^M \big( \dim \xxF_{\sigma}^{\geq n}  - \dim \xxF_{\sigma}^{\geq n + 1}\big)\\
&= \dim \xxF_{\sigma}^{\geq -N}  - \dim \xxF_{\sigma}^{\geq M + 1} = \dim \xxF_{\sigma}.
\end{align*}
Hence $\overline \xxF$ has the same Hilbert function as $\xxF$, so that the sum of the graded pieces $[q_n\colon \xxH^n \twoheadrightarrow \xxF^{[n]}]$ yields a point $[\overline{q} = \oplus_{n} q_n \colon \xxH \twoheadrightarrow \overline{\xxF}]  \in \xxGQuotH$.
It has the property that it is the limit of the action of $\lambda(t)$ on $[q]$ when $t$ tends to infinity:

\begin{lem}  \label{limq}
Let $[q \colon \xxH \twoheadrightarrow\xxF]  \in \xxGQuotH$, let $\lambda: \kxC^* \to \Gamma$ be a $1$--parameter subgroup, and let $[\overline{q}]$ be the quotient of $\xxGQuotH$ defined as above. Then \linebreak
$[\overline{q}] = \lim_{t \to 0} [q]\cdot\lambda(t)^{-1} = \lim_{t \to \infty} [q]\cdot\lambda(t)$. 
\end{lem}
\begin{proof}
The proof works analogously to \cite[Lemma 4.4.3]{HL:2010}. The main difference is a minus sign, which occurs since we consider descending filtrations while \cite{HL:2010} work with ascending filtrations. Therefore, we obtain the limit at infinity instead of zero. Consult \cite[Lemma 3.3.2]{Bec:thesis} for the details.
\end{proof}

The description of $[\overline{q}]$ as a limit of $[q]\cdot\lambda(t)$ yields that it is a fixed point of the action of $\lambda$. Hence there is an action of $\lambda$ on the fibre 
$$\xLl_{\chi}([\overline{q}]) = \bigotimes_{\sigma \in D} \det (\overline{\xxF}_{\sigma})^{\kappa_{\sigma}} = \bigotimes_{\sigma \in D} \det \big(\bigoplus_{n \in \intZ} \xxF^{[n]}_{\sigma}\big)^{\kappa_{\sigma}} = \bigotimes_{\sigma \in D} \bigotimes_{n \in \intZ} \det ( \xxF^{[n]}_{\sigma})^{\kappa_{\sigma}}.$$
We examine this action in the sequel so as to gain some criteria for the \linebreak GIT--(semi)stability of $[q]$.

\subsection{GIT--(semi)stability} \label{GITstab1}

\subsubsection{A numerical criterion}  \label{GITstab2}
Let us start by recalling the notion of (semi)stability in the GIT--sense as defined in \cite[Definition 1.7]{GIT:1994}.

\begin{definition}
Let $\Gamma$ be a reductive algebraic group, let $S$ be a $\Gamma$--scheme, and let $\xLl$ be a $\Gamma$--linearized line bundle on $S$. Then:
\begin{itemize}
\item A (closed) point $s \in S$ is semistable if there exists a section $\sigma \in H^0(X,\xLl^n)^\Gamma$ for some $n$, such that $\sigma(s) \neq 0$, and $X_\sigma:=\{ s \in S | \sigma(s) \neq 0\}$ is affine.  
\item A (closed) point $s \in S$ is stable if it is semistable, $\Gamma.s$ is closed in $X_\sigma$, and the isotropy group $\Gamma_s$ is finite. 
\end{itemize}
\end{definition} 

Note that the definition of GIT--stable point given here differs slightly from the one in \cite{GIT:1994} since Mumford does not require for $\Gamma_s$ to be finite. It is obvious from the definition that the set of GIT--(semi)stable points is open and $\Gamma$--invariant.

We have seen in Section \ref{GITparameters} that the invariant Quot scheme $\xxGQuotH$ is equipped with a $\Gamma$--linearized ample line bundle $\xLl_\chi$. 
Hence, by \cite[Amplification 1.8]{GIT:1994}, there exists some integer $N \geq 1$, a finite dimensional $\Gamma$--submodule $W \subset H^0(\xxGQuotH,\xLl_\chi^N)$, and a $\Gamma$--equivariant locally closed immersion 
$$\iota : \xxGQuotH \hookrightarrow \prP(W)$$ 
such that $\iota^* \ssh_{\prP(W)}(1) \cong \xLl_\chi^N$ and, if we identify $\xxGQuotH$ with a subscheme of $\prP(W)$ via $\iota$, then  
\begin{equation}  \label{Ampli1.8}
\xxGQuotH^{(s)s}(\xLl_\chi)=\xxGQuotH \cap \prP(W)^{(s)s}(\ssh_{\prP(W)}(1)),
\end{equation}
where we denote by $\xxGQuotH^{(s)s}(\xLl_\chi)$ and by $\prP(W)^{(s)s}(\ssh_{\prP(W)}(1))$ the open subset of GIT--(semi)stable points with respect to the linearized line bundle $\xLl\chi$ and $\ssh_{\prP(W)}(1)$ respectively. In fact, \cite[Amplification 1.8]{GIT:1994} gives equality \eqref{Ampli1.8} only for GIT--stable points, but the same proof gives the result also for GIT--semistable points.

As $\prP(W)$ is a proper scheme, we can use Mumford's numerical criterion to determine the GIT--(semi)stable locus of $\xxGQuotH$ with respect to $\xLl_\chi$. We first adapt Mumford's definition \cite[Definition 2.2]{GIT:1994} to our situation:
  
\begin{definition}
For $[q \colon \xxH \twoheadrightarrow \xxF] \in \xxGQuotH$ and any $1$--parameter subgroup $\lambda$ of $\Gamma$ we define $\mu_{\xLl_{\chi}}(q,\lambda)$ as the weight of $\lambda$ on $\xLl_{\chi}([\overline{q}])$.
\end{definition}

It follows from the discussion above that, in our situation, Mumford's numerical criterion \cite[Theorem 2.1]{GIT:1994} can be formulated as follows:

\begin{prop}[Mumford's numerical criterion]\hfill\par
\noindent The point $[q \colon \xxH \twoheadrightarrow \xxF] \in \xxGQuotH$ is GIT--(semi)stable with respect to the twisted line bundle $\xLl_{\chi}$ if and only if for every non--trivial $1$--parameter subgroup $\lambda\colon \kxC^* \to \Gamma$ we have $\mu_{\xLl_{\chi}}(q,\lambda) \sgeq 0$.
\end{prop}

%\begin{remark}
%In the case of vector bundles, GIT--(semi)stability is equivalent to the condition $\mu(q,\lambda)\sleq 0$ \cite[Theorem 4.2.11]{HL:2010} when $\mu$ is defined via the weight of $\lambda$ on the fibre of $\xLl$ at the limit at zero. As we consider the limit at infinity, or equivalently the limit of the inverse $1$--parameter subgroup at zero, we have GIT--(semi)stability exactly when the negative weight is $\sleq 0$, i.e.~$\mu_{\xLl_{\chi}}(q,\lambda)\sgeq 0$.
%\end{remark}

Now we establish some expressions for $\mu_{\xLl_{\chi}}(q,\lambda)$ in terms of $\kappa$ and $\chi$:

\begin{lem}
The weight of the action of $\kxC^*$ via $\lambda$ on $\xLl_{\chi}[\overline{q}]$ is
\begin{align*}
\mu_{\xLl_{\chi}}(\overline{q},\lambda) &= \sum_{n \in \intZ}n\biggl(\sum_{\sigma \in D}\kappa_{\sigma}\cdot \dim_{\kxC}(\xxF^{[n]}_{\sigma}) + \sum_{\rho \in D_-} \chi_{\rho} \cdot \dim_{\kxC}(A_{\rho}^n)\biggr) \\
&=: \sum_{n \in \intZ}n\bigl(\kappa(\xxF^{[n]}) + \chi(A^n)\big).
\end{align*}
\end{lem}
\begin{proof}
The weight $\mu_{\xLl_{\chi}}(\overline{q},\lambda)$ is the exponent in the identity
$$\lambda(t)|_{\xLl_{\chi}([\overline{q}])} = \lambda(t)|_{\bigotimes\limits_{\sigma \in D} \bigotimes\limits_{n \in \intZ} \det ( \xxF^{[n]}_{\sigma})^{\kappa_{\sigma}}} = t^{\mu_{\xLl_{\chi}}(\overline{q},\lambda)}\cdot\xid_{\xLl_{\chi}([\overline{q}])}.$$
This number splits into a sum $\mu_{\xLl_{\chi}}(\overline{q},\lambda) = m + m_{\chi}$, where $m$ is the weight on the fibre of the original line bundle $\xLl([\overline{q}])$ and $m_{\chi}$ comes from the twist with the character $\chi$.

Since the weight of $\lambda$ on $\xxF^{[n]}_{\sigma}$ is $n$, for its weight on the determinant $\det(\xxF^{[n]}_{\sigma})^{\kappa_{\sigma}}$ we obtain $n\cdot \dim(\xxF^{[n]}_{\sigma}) \cdot \kappa_{\sigma}$. The weights on the factors of the tensor products over $D$ and $\intZ$ translate to a sum of the weights, thus
$m = \sum_{\sigma \in D} \sum_{n \in \intZ}n\cdot\kappa_{\sigma}\cdot\dim \xxF^{[n]}_{\sigma}$.

The  $\lambda(t)_{\rho}$ are diagonal matrices of size $(\dim A_{\rho}) \times (\dim A_{\rho})$ with entries $t^n$ according to the decomposition $A_{\rho} = \bigoplus_{n \in \intZ} A_{\rho}^n$. The twist by the character $\chi$ is given by taking the product of the determinants of the $\lambda(t)_{\rho}$ to the $\chi_{\rho}$'s power. Thus we have 
$$t^{m_{\chi}} = \prod_{\rho \in D_-}\det(\lambda(t)_{\rho})^{\chi_{\rho}} = \prod_{\rho \in D_-} \prod_{n \in \intZ} t^{n\cdot \dim(A_{\rho}^n) \cdot \chi_{\rho}},$$
and $m_{\chi} = \sum_{\rho \in D_-} \sum_{n \in \intZ} n\cdot\chi_{\rho} \cdot \dim(A_{\rho}^n)$.
%
%Together, this yields
%$$\mu_{\xLl_{\chi}}(\overline{q},\lambda) = \sum_{n \in \intZ}n \biggl( \sum_{\sigma \in D}\kappa_{\sigma}\cdot\dim \xxF^{[n]}_{\sigma} + \sum_{\rho \in D_-}\chi_{\rho}\cdot \dim A_{\rho}^n \biggr).$$
%The weight of $\lambda$ on $\xxF^{[n]}_{\sigma}$ is $n$, so
%\begin{align*}
%\lambda(t)|_{\xLl_{\chi}([\overline{q}])} &= \lambda(t)|_{\bigotimes\limits_{\sigma \in D} \bigotimes\limits_{n \in \intZ} \det ( \xxF^{[n]}_{\sigma})^{\kappa_{\sigma}}} = \prod_{\sigma \in D} \prod_{n \in \intZ} \det(t^n \cdot \xid_{\xxF^{[n]}_{\sigma}})^{\kappa_{\sigma}} \cdot \prod_{\rho \in D_-}\det(\lambda(t)_{\rho})^{\chi_{\rho}}\\
%&= \prod_{\sigma \in D} \prod_{n \in \intZ} t^{n\cdot\dim \xxF^{[n]}_{\sigma} \cdot \kappa_{\sigma}} \det (\xid_{\xxF^{[n]}_{\sigma}})^{\kappa_{\sigma}} \cdot \prod_{\rho \in D_-}t^{\sum_{n \in \intZ}n\cdot \dim A_{\rho}^n\cdot\chi_{\rho}}\\
%&= t^{\sum_{n \in \intZ} \sum_{\sigma \in D} n\cdot\kappa_{\sigma}\cdot\dim \xxF^{[n]}_{\sigma}} \prod_{\sigma \in D} \prod_{n \in \intZ}  \xid_{(\det\xxF^{[n]}_{\sigma})^{\kappa_{\sigma}} } \cdot t^{\sum_{n \in \intZ}\sum_{\rho \in D_-}n\cdot\chi_{\rho}\cdot \dim A_{\rho}^n}\\
%&= t^{\sum_{n \in \intZ}n \left( \sum_{\sigma \in D}\kappa_{\sigma}\cdot\dim \xxF^{[n]}_{\sigma} + \sum_{\rho \in D_-}\chi_{\rho}\cdot \dim A_{\rho}^n \right)} \cdot \xid_{\xLl_{\chi}([\overline{q}])}.
%\end{align*}
\end{proof}

%Let us note that for every $n,m \in \intZ$, we have $\kappa(\xxF^{[n]} \oplus \xxF^{[m]})=\kappa(\xxF^{[n]})+\kappa(\xxF^{[m]})$, resp. $\chi(A^n \oplus A^m)=\chi(A^n)+\chi(A^m)$. 
Generalizing the calculation before Proposition 3.1 in \cite{Kin:1994}, we obtain another formula for $\mu_{\xLl_{\chi}}(\overline{q},\lambda)$:  
\begin{prop}
In terms of the filtration corresponding to a $1$--parameter subgroup $\lambda$, we have
$$\mu_{\xLl_{\chi}}(\overline{q},\lambda) = \sum_{n=-N+1}^M \left( \kappa(\xxF^{\geq n}) + \chi(A^{\geq n}) \right) - N \cdot \kappa(\xxF),$$
where $-N$ is the minimal and $M$ the maximal occurring weight.
\end{prop}

\begin{proof}
Using the fact that $\xxF^{\geq n} = 0$, $A^{\geq n} = 0$ for $n > M$ and $\xxF^{\geq n} = \xxF$, $A^{\geq n} = A$ for $n \leq -N$, this is an easy calculation. It is accomplished in \cite[Proposition 3.4.4]{Bec:thesis}.%so we can use Lemma \ref{hilfsrechnung} twice, setting $B = \xxF$, $\varphi = \kappa$ and $B = A$, $\varphi = \chi$, respectively. This yields
%\begin{align*}
%\sum_{n \in \intZ}n\bigl(\kappa(\xxF^{[n]}) + \chi(A^n)\big)
%&\overset{\eqref{gradfiltid}}{=} \sum_{n = -N+1}^{M} \left( \kappa(\xxF^{\geq n}) + \chi(A^{\geq n}) \right) -  N \cdot \big( \kappa(\xxF) + \underbrace{\chi(A)}_{=0} \big)\\
%&= \sum_{n = -N+1}^{M} \left( \kappa(\xxF^{\geq n}) + \chi(A^{\geq n}) \right) -  N \cdot \kappa(\xxF).
%\end{align*}
\end{proof}

\subsubsection{$1$--step filtrations}

Next we analyse the stability condition for $1$--step filtrations in order to simplify the condition for GIT--(semi)stability:

Let $\lambda$ be a $1$--parameter subgroup such that the corresponding filtration is a one-step filtration $A \supsetneq A' \supsetneq 0$, and let $A''$ be a complement of $A'$ in $A$. Denote by $n'$ and $n''$ the weights of $A'$ and $A''$ for the action of $\lambda$ on $A$, respectively. 

As $n' \dim A'+n'' \dim A''=0$ by equation \eqref{sumAm} and $n''<n'$, we have $n'=\dim A''=\dim A- \dim A'$ and $n''=-\dim A'$ up to a multiple in $\frac{1}{gcd(\dim A', \dim A)} \natN$. Moreover:

\begin{align*}
& A^{n'} = A',\\
& A^{n''} = A'' \cong A/A',\\
& \xxF^{[n']} = q\big(\bigoplus_{\rho \in D_-} A'_{\rho} \tensor V_{\rho} \tensor \ssh_X\big) =: \xxF',\\
& \xxF^{[n'']} = q(\bigoplus_{\rho \in D_-} (A'_{\rho} \oplus A''_{\rho}) \tensor V_{\rho} \tensor \ssh_X)/\xxF' = \xxF/\xxF'.
\end{align*}
As $\chi(A)=\sum_{\rho \in D_-} \chi_\rho h(\rho)=0$ by assumption on $\chi$, we get
\begin{align*}
\mu_{\xLl_{\chi}}(q,\lambda) &= n'\cdot \bigl( \kappa(\xxF') + \chi(A') \bigr) + n'' \cdot \bigl( \underbrace{\kappa(\xxF/\xxF')}_{\kappa(\xxF)-\kappa(\xxF')} + \;\underbrace{\chi(A/A')}_{-\chi(A')} \;\bigr)\\
&= (n'-n'') \cdot \bigl( \kappa(\xxF') + \chi(A') \bigr) + n'' \cdot \kappa(\xxF).
\end{align*}
Thus we obtain the following criterion for $\mu_{\xLl_{\chi}}(q, \lambda)$ to be positive:
\begin{equation}\label{onestepstab}
\mu_{\xLl_{\chi}}(q,\lambda) \sgeq 0 \quad \Longleftrightarrow \quad \mu(A') := \dim A \cdot \bigl( \kappa(\xxF') + \chi(A') \bigr) - \dim A' \cdot \kappa(\xxF) \sgeq 0. 
\end{equation}

The next proposition gives a criterion for GIT--(semi)stability in terms of graded subspaces:

\begin{prop}  \label{ineqmuA}
An element $[q\colon \xxH \twoheadrightarrow \xxF] \in \xxGQuotH$ is GIT--(semi)stable if and only if for every graded subspace $0 \neq A' \subsetneq A$, that is $A'=\bigoplus_{\rho \in D_-} A'_\rho$ with $A'_\rho \subset A_\rho$, and for every $\xxF' := q\big(\bigoplus_{\rho \in D_-} A'_{\rho} \tensor V_{\rho} \tensor \ssh_X\big)$ the inequality $\mu(A') := \dim A \cdot \bigl( \kappa(\xxF') + \chi(A') \bigr) - \dim A' \cdot \kappa(\xxF) \sgeq 0$ holds.
\end{prop}

\begin{proof}
``$\Rightarrow$'': Considering a $1$--parameter subgroup corresponding to the $1$--step filtration 
$$0 \neq A' \subsetneq A,$$ 
this follows from \eqref{onestepstab} and Mumford's numerical criterion.

``$\Leftarrow$'': Let $\lambda$ be any non--trivial $1$--parameter subgroup. By Mumford's numerical criterion we have to show that $\mu_{\xLl_{\chi}}(q,\lambda) \sgeq 0$. Let $-N$ denote the minimal and $M$ the maximal occurring weight. For every $n \in \{-N+1,\hdots,M\}$ let us consider the graded subspace $A \supsetneq A^{\geq n} \supsetneq 0$. By assumption $\kappa(\xxF^{\geq n}) + \chi(A^{\geq n}) > \frac{\dim A^{\geq n}}{\dim A}\cdot\kappa(\xxF)$. 
This yields 
\begin{align*}
\mu_{\xLl_{\chi}}(q,\lambda) &= \sum_{n = -N+1}^M \bigl( \kappa(\xxF^{\geq n}) + \chi(A^{\geq n}) \bigr) - N\cdot \kappa(\xxF) \\
& \sgeq \sum_{n = -N+1}^M \dim A^{\geq n} \cdot \frac{\kappa(\xxF)}{\dim A} - N\cdot \kappa(\xxF) \quad = \quad 0,%\\
%&= N \cdot \dim A \cdot \frac{\kappa(\xxF)}{\dim A} - N\cdot \kappa(\xxF) = 0,
\end{align*}
since one checks that %by Lemma \ref{hilfsrechnung} with $B = A$, $\varphi \equiv 1$ we have
$$\sum_{n = -N+1}^M \dim A^{\geq n} = \underbrace{\sum_{n \in \intZ} n \cdot \dim A^n}_{= 0 \text{ by } \eqref{sumAm}} + \; N \cdot \dim A = N \cdot \dim A.$$
%the calculation
%\begin{align*}
%\sum_{n \in \intZ} n \cdot \dim A^n &= \sum_{n \in \intZ} n \cdot (\dim A^{\geq n} - \dim A^{\geq n+1})\\
%&= \sum_{n = -N}^M n \cdot \dim A^{\geq n} - \sum_{n = -N+1}^{M+1} (n-1) \cdot \dim A^{\geq n}\\
%&=\sum_{n = -N+1}^M (n-(n-1)) \cdot \dim A^{\geq n} + (-N) \cdot \underbrace{\dim A^{\geq -N}}_{= \dim A} - (M + 1 - 1) \cdot \underbrace{\dim A^{\geq M}}_{= 0}\\
%&= \sum_{n = -N+1}^M \dim A^{\geq n} - N \cdot \dim A.
%\end{align*}

This shows that $[q]$ is GIT--(semi)stable.
\end{proof}

% *** 4 ***************************************************************************************************************************************
\section{The connection between the stability conditions}\label{chcorres}

As we want to construct the moduli space of $\theta$--stable $(G,h)$--constellations on an affine $G$--scheme $X$ as an open subset of the GIT--quotient $\xxGQuotH^{ss} \red_{\xLl_{\chi}} \Gamma$, first of all we determine the elements in $\xxGQuotH$ originating from $(G,h)$--con\-stellations in Subsection \ref{quotconst}. It turns out that every GIT--semistable quotient can indeed be obtained from a $(G,h)$--constellation in a particular way, so that we can define a functor $\fsMck$ of flat families of GIT--stable $(G,h)$--constellations. We compare $\fsMck$ to the \mbox{functor} $\fsMth$ of flat families of $\theta$--stable $(G,h)$--constel\-lations. Therefore, in Subsection \ref{FAcorres} we establish a correspondence between the $G$--equi\-variant coherent subsheaves generated in $D_-$ of a $(G,h)$--constellation $\xxF$ and the saturated graded subspaces of $A = \bigoplus_{\rho \in D_-}A_{\rho}$ defining subsheaves of $\xxH$. This leads us to the definition of a new stability condition $\widetilde \theta$ on $(G,h)$--constellations which coincides with GIT--stability for $(G,h)$--constellations generated in $D_-$. This reduces our examination of the stability conditions to a comparison of $\theta$ and $\widetilde \theta$, which look very similar for a certain choice of the GIT--parameters $\chi$ and $\kappa$. Indeed, in Subsection \ref{comparison} we show that $\theta$ is a limit of the $\widetilde \theta$, when the finite subset $D \subset \Irr G$ in the definition of $\widetilde \theta$ varies. Furthermore, we find out that $\theta$--stability implies $\widetilde \theta$--stability and hence GIT--stability, so that the functor of $\theta$--stable $(G,h)$--constellations is a subfunctor of the functor of GIT--stable $(G,h)$--constellations.

\subsection{Quotients originating from $(G,h)$--constellations} \label{quotconst}

To determine those points in the invariant Quot scheme which originate from $\theta$--semistable $(G,h)$--con\-stellations, we analyse the quotient map for these elements first.

Let $\xxF$ be a $\theta$--semistable $(G,h)$--constellation generated in $D_-$. By Theorem \ref{finiteness}, this is the case for instance if $\xxF$ is $\theta$--stable or if $D_0=\{ \rho \in \Irr G \ |\ h(\rho)=0\}$. By Section \ref{sfiniteness} we can write $\xxF$ as a quotient of
$$\xxH := \bigoplus_{\rho \in D_-} A_{\rho} \tensor V_{\rho} \tensor \ssh_X,$$
where $A_{\rho} = \kxC^{h(\rho)}$ and $D_-$ is the finite subset of $\Irr G$ where $\theta$ takes negative values.
Since $\xxF_{\rho} = \mathcal{H}om_G(V_{\rho}, H^0(\xxF))$ we have natural evaluation maps
$$ev_{\rho}\colon \xxF_{\rho} \tensor V_{\rho} \tensor \ssh_X \to \xxF,\ \;\alpha \otimes v \otimes f \mapsto f \cdot \alpha(v)$$
and $\xxF$ is generated as an $\ssh_X$--module by the images of $ev_{\rho}$, $\rho \in D_-$ by assumption. Choosing a basis of each $\xxF_{\rho}$, i.e.~fixing an isomorphism $\psi_{\rho}\colon A_{\rho} \to \xxF_{\rho}$, and composing it with the evaluation map, we obtain
\begin{equation}\label{quotofconst} q_{\rho}\colon A_{\rho} \tensor V_{\rho} \tensor \ssh_X \to \xxF,\ \;a \otimes v \otimes f \mapsto f \cdot \psi_{\rho}(a)(v).\end{equation}
Their sum
$$q := \underset{\rho \in D_-}{\oplus} q_{\rho}\colon \xxH = \bigoplus_{\rho \in D_-} A_{\rho} \tensor V_{\rho} \tensor \ssh_X \to \xxF$$
gives us a point $[q\colon \xxH \twoheadrightarrow \xxF] \in \xxGQuotH$ with the property that the map
\begin{equation} \label{iso} \varphi_{\rho}\colon A_{\rho} \to \xxF_{\rho} = \mathcal{H}om_G(V_{\rho},H^0(\xxF)),\ \; a \mapsto (v \mapsto q(a \otimes v \otimes 1)), \end{equation} 
is just the isomorphism $\psi_{\rho}$ since, for $a \in A_{\rho}$ and $v \in V_{\rho}$, we have
$$\varphi_{\rho}(a)(v) = q(a \otimes v \otimes 1) = 1 \cdot \psi_{\rho}(a)(v) = \psi_{\rho}(a)(v).$$

The point $[q\colon \xxH \twoheadrightarrow \xxF] \in \xxGQuotH$ constructed this way depends on the choice of the isomorphisms $\psi_{\rho}$. Any other choice differs from $\psi_{\rho}$ by an element in $Gl({A_{\rho}})$, so that a $(G,h)$--constellation can be seen as an element in the quotient of $\xxGQuotH$ by the gauge group $\Gamma$ defined by \eqref{def Gamma}. We will make this more precise in Section \ref{modspace}.

Conversely, for any element $[q\colon \xxH \twoheadrightarrow \xxF] \in \xxGQuotH$, the quotient $\xxF$ is a $G$--equivariant coherent $\ssh_X$--module with isotypic decomposition isomorphic to $R_h$, so it is a $(G,h)$--constellation. However, the induced maps $\varphi_{\rho}$ need not be isomorphisms so that $[q]$ need not originate from a $(G,h)$--constellation as above even if $\xxF$ is $\theta$--stable. Since we want to determine a moduli space $\sMth$ of $\theta$--stable $(G,h)$--constellations as a subscheme of $\xxGQuotH^{ss}\red_{\xLl_{\chi}}\Gamma$, we are interested in exploring which quotient maps $q$ do indeed arise from a $(G,h)$--constellation.

\begin{lem} \label{GITssiso}
Let $[q \colon \xxH \twoheadrightarrow \xxF] \in \xxGQuotH$ be GIT--semistable, and let $\rho \in D_-$. If $\chi_{\rho} < \frac{\kappa(\xxF)}{\dim A}$, then $\varphi_{\rho}\colon A_{\rho} \to \xxF_{\rho}$ is an isomorphism. 
\end{lem}

\begin{proof}
Fix $\rho \in D_-$ and let $K_{\rho} := \ker \varphi_{\rho}$. If $\varphi_{\rho}$ is not injective, then one can find a $1$--parameter subgroup $\lambda$ such that $A \supset K_{\rho} \supsetneq 0$ is a ($1$--step) filtration corresponding to $\lambda$. For the induced subsheaf we obtain $\xxF' = q(K_{\rho} \tensor V_{\rho} \tensor \ssh_X) =\ssh_X \cdot \varphi_{\rho}(K_{\rho})(V_\rho) = 0$, so that
\begin{align*}
\mu(K_{\rho}) &= \dim A \cdot \bigl( \kappa(0) + \chi(K_{\rho}) \bigr) - \dim K_{\rho} \cdot \kappa(\xxF) \\
&=\dim A \cdot \chi_{\rho}\dim K_{\rho} - \dim K_{\rho} \cdot \kappa(\xxF)\\
&=\dim K_{\rho} \cdot (\dim A \cdot \chi_{\rho} - \kappa(\xxF) ) \;<\; 0
\end{align*}
by the assumption on $\chi_{\rho}$.

This is a contradiction to semistability by Proposition \ref{ineqmuA}, so $\ker \varphi_{\rho}$ has to be $0$. As $A_{\rho}$ and $\xxF_{\rho}$ have the same dimension $h(\rho)$, this implies that $\varphi_{\rho}$ is an isomorphism.
\end{proof}

Lemma \ref{GITssiso} means that, if $\chi_{\rho} < \frac{\kappa(\xxF)}{\dim A}$, then for every GIT--semistable quotient $[q\colon \xxH \twoheadrightarrow \xxF] \in \xxGQuotH$ the $q_{\rho},\ \rho \in D_-,$ are of the form \eqref{quotofconst}. In this sense, $[q]$ arises from a $(G,h)$--constellation.

If for a $(G,h)$--constellation $\xxF$ and a choice of isomorphisms $(\psi_{\rho})_{\rho \in D_-}$ the corresponding point is GIT--(semi)stable, then the same is true for any other choice of isomorphisms by $\Gamma$--invariance of GIT--(semi)stable points. 
Thus it makes sense to deal with GIT--(semi)stable $(G,h)$--constellations:

\begin{definition}  \label{GITfunctors}
A $(G,h)$--constellation $\xxF$ is \textit{GIT--(semi)stable}, if for some (and hence any) choice of isomorphisms $(\psi_{\rho})_{\rho \in D_-}$ the corresponding point as defined in \eqref{quotofconst} is GIT--(semi)stable.
Let
\begin{align*}
&\fMck\colon \text{(Sch/$\kxC$)}^{\text{op}} \to \text{(Set)} \\
&S \mapsto \{\xxFf \text{ an } S\text{--flat family of GIT--semistable $(G,h)$--constellations on $X \times S$}\}/_{\cong} \\
&(f\colon S' \to S) \mapsto \big(\fMck(S) \to \fMck(S'), \xxFf \mapsto (\xid_X\times f)^*\xxFf\big),\\
\intertext{and}
&\fsMck\colon \text{(Sch/$\kxC$)}^{\text{op}} \to \text{(Set)} \\
&S \mapsto \{\xxFf \text{ an } S\text{--flat family of GIT--stable $(G,h)$--constellations on $X \times S$}\}/_{\cong} \\
&(f\colon S' \to S) \mapsto \big(\fsMck(S) \to \fsMck(S'), \xxFf \mapsto (\xid_X\times f)^*\xxFf\big)
\end{align*}
be the moduli functors of GIT--semistable and GIT--stable $(G,h)$--constellations on $X$ generated in $D_-$, respectively.%, and further 
\end{definition}

From the discussion above we expect that the quotients $\xxGQuotH^{ss}\red_{\xLl_{\chi}}\Gamma$ and $\xxGQuotH^s/\Gamma$ corepresent these functors. We will prove this in Section \ref{corepr}.

%**** 4.2 *************************************************************************************************************************************************
\subsection[Correspondence between $A' \subset A$ and $\xxF' \subset \xxF$]{Correspondence between saturated graded subspaces of $A$ and $G$--equivariant subsheaves of $\xxF$ generated in $D_-$}\label{FAcorres}

If the map $A_{\rho} \to \xxF_{\rho}$ is injective and hence an isomorphism, we may establish a correspondence between subsheaves of the $(G,h)$--constellation $\xxF$ generated in $D_-$ and saturated graded subspaces of $A$. By Lemma \ref{GITssiso} this correspondence applies to GIT--semistable elements.
First we begin with some graded subspace $A' \subset A$, i.e.~we have subspaces $A_{\rho}' \subset A_{\rho}$ for every $\rho \in D_-$. Let 
\begin{equation}\label{FausA}
\xxF' := q \big( \bigoplus_{\rho \in D_-} A_{\rho}' \tensor V_{\rho} \tensor \ssh_X \big) =  \ssh_X \cdot \big( \sum_{\rho \in D_-} \varphi_{\rho}(A_{\rho}')(V_\rho)\big)
\end{equation}
be the $G$--equivariant subsheaf of $\xxF$ generated by the $\varphi_{\rho}(A_{\rho}')$, $\rho \in D_-$. Since $\varphi_{\rho}|_{A_{\rho}'}$ is injective we have $\dim A_{\rho}' \leq \dim \xxF_{\rho}'$ for every $\rho \in D_-$.
Further, we define
$$\widetilde A_{\rho}' := \varphi_{\rho}^{-1}(\xxF_{\rho}'), \qquad \widetilde A' := \bigoplus_{\rho \in D_-} \widetilde A_{\rho}'.$$
Then we have
\begin{itemize}\renewcommand{\labelitemi}{\textbullet}
 \item $\dim \widetilde A_{\rho}' = \dim \xxF_{\rho}' =: h'(\rho)$ since $\varphi_{\rho}$ is an isomorphism,
 \item $\widetilde A_{\rho}' = \varphi_{\rho}^{-1}\big(\xxF'_{\rho}\big) \supset \varphi_{\rho}^{-1}\big(\varphi_{\rho}(A_{\rho}')\big) = A_{\rho}'$,
 \item $q\big( \bigoplus_{\rho \in D_-} \widetilde A_{\rho}' \tensor V_{\rho} \tensor \ssh_X\big) = \ssh_X \cdot \big( \sum_{\rho \in D_-} \varphi_{\rho}(\widetilde A_{\rho}')(V_\rho) \big) = \xxF'$, since $\varphi_{\rho}(\widetilde A_{\rho}') = \xxF_{\rho}'$ if $\rho \in D_-$ and $\xxF'$ is generated in $D_-$.
\end{itemize}
For this reason, $\widetilde A'$ is called the saturation of $A'$.

Inspired by this correspondence we define a new function, which describes \linebreak GIT--(semi)stability in terms of the $\xxF'$ instead of the $A'$:
\begin{definition} \label{thetatilde} Let $\xxF$ be any $(G,h)$--constellation, $\xxF' \subset \xxF$ a $G$--equivariant coherent subsheaf, $h'(\rho) := \dim \xxF_{\rho}'$. Let $\widetilde \theta\colon \Coh^G(X) \to \ratQ$ be the function
$$ \widetilde \theta(\xxF') := \sum_{\rho \in D_-}\Big(\kappa_{\rho} + \chi_{\rho} -\frac{\kappa(\xxF)}{\dim A}\Big)h'(\rho)\;\; + \sum_{\sigma \in D\setminus D_-}\kappa_{\sigma}h'(\sigma),$$
where we recall that $\kappa(\xxF)=\sum_{\sigma \in D} \kappa_\sigma h(\sigma)$ and $A=\bigoplus_{\rho \in D_-} A_\rho$.
\end{definition}

In the above setting, if $\xxF'$ is generated in $D_-$, then we have $h'(\rho) = \dim \widetilde A'_{\rho}$. Comparing this definition to expression \eqref{onestepstab} we find
\begin{equation} \label{thetatildemu} \dim A \cdot \widetilde \theta(\xxF') 
= \mu(\widetilde A').\end{equation}

\begin{remark}
Since the notion of GIT--stability on $\xxGQuotH$ depends on the embedding into a product of Grassmannians, the definition of $\widetilde \theta$ depends on the choice of the finite subset $D \subset \Irr G$. If there is any ambiguity about $D$ we write $\widetilde \theta_D$ instead of $\widetilde \theta$.
\end{remark}

The next theorem reduces the examination of the relation between $\theta$--(semi)stabi\-lity and GIT--(semi)stability to the comparison of $\theta$ and $\widetilde \theta$ for sheaves generated in $D_-$.

\begin{thm} \label{equivGITstab}
Let $\xxF$ be a $(G,h)$--constellation generated in $D_-$, and suppose that for every $\rho \in D_-$ we have $\chi_{\rho} < \frac{\kappa(\xxF)}{\dim A}$.
Then $\xxF$ is GIT--(semi)stable if and only if $\xxF$ is $\widetilde\theta$--(semi)stable.
\end{thm}
\begin{proof}
``$\Rightarrow$'': 
Let $\xxF'$ be a $G$--equivariant subsheaf of $\xxF$ generated in $D_-$. As above we consider $\widetilde A' = \bigoplus_{\rho \in D_-}\varphi_{\rho}^{-1}(\xxF_{\rho}')$. Then we have $\widetilde \theta(\xxF') = \frac{\mu(\widetilde A')}{\dim A} \sgeq 0$ by GIT--(semi)stability. 
% and the subsheaf $\widetilde \xxF'$ of $\xxF'$ generated by the $\xxF'_{\rho}$, $\rho \in D_-$. Thus $\widetilde h'(\rho) := \dim \widetilde\xxF_{\rho}' = h'(\rho)$ for $\rho \in D_-$ and $\widetilde h'(\rho) \leq h'(\rho)$ for $\rho \in \Irr G \setminus D_-$. As $\widetilde \xxF'$ is generated in $D_-$, we have $\widetilde \theta(\widetilde \xxF') = \frac{\mu(\widetilde A')}{\dim A} \sgeq 0$ by GIT--(semi)stability. 
%For $\xxF'$ this yields
%\begin{align*}
%\widetilde\theta(\xxF') &= \sum_{\rho \in D_-}\Big(\kappa_{\rho} + \chi_{\rho} -\frac{\kappa(\xxF)}{\dim A}\Big)\underbrace{h'(\rho)}_{= \widetilde h'(\rho)} \; + \sum_{\sigma \in D \setminus D_-}\underbrace{\kappa_{\sigma}}_{>0}\underbrace{h'(\sigma)}_{\geq \widetilde h'(\sigma)} \\
%&\geq \sum_{\rho \in D_-}\Big(\kappa_{\rho} + \chi_{\rho} -\frac{\kappa(\xxF)}{\dim A}\Big)\, \widetilde h'(\rho) \; + \sum_{\sigma \in D \setminus D_-}\kappa_{\sigma}\widetilde h'(\sigma) = \widetilde\theta(\widetilde\xxF') \;\sgeq\; 0.
%\end{align*}

``$\Leftarrow$'': Let $A' \subset A$ be a graded subspace. As in \eqref{FausA} we construct $\xxF'$ and $\widetilde A' \supset A'$. By $\widetilde\theta$--(semi)stability we have 
$\mu(\widetilde A') = \dim A \cdot \widetilde \theta(\xxF') \sgeq 0$. If $\widetilde A'=A'$, then equality \eqref{thetatildemu} gives the result. Otherwise, we obtain
\begin{align*}
\chi(\widetilde A') - \chi(A') = \chi(\widetilde A'/A') &= \sum_{\rho \in D_-} \chi_{\rho} \cdot \dim(\widetilde A'/A')_{\rho}\\
&< \sum_{\rho \in D_-} \frac{\kappa(\xxF)}{\dim A} \cdot \dim(\widetilde A'/A')_{\rho}\\
&= \frac{\kappa(\xxF)\cdot \dim(\widetilde A'/A')}{\dim A} = \frac{\dim\widetilde A'-\dim A'}{\dim A}\cdot\kappa(\xxF).
\end{align*}
Separating $\widetilde A'$ and $A'$ and multiplying by $\dim A$ yields
$$\dim A \cdot \chi(\widetilde A') - \dim \widetilde A' \cdot \kappa(\xxF) < \dim A \cdot \chi(A') - \dim A' \cdot \kappa(\xxF),$$
so that
\begin{align*}
\mu(A') &= \dim A \cdot (\kappa(\xxF') + \chi(A')) - \dim A' \cdot \kappa(\xxF) \\
&> \dim A \cdot (\kappa(\xxF') + \chi(\widetilde A')) - \dim \widetilde A' \cdot \kappa(\xxF) = \mu(\widetilde A') \geq 0.
\end{align*}
\end{proof}

If we could show that
\begin{equation}\label{thequivtth}\widetilde \theta (\xxF') \sgeq 0 \quad \Longleftrightarrow \quad \theta (\xxF') \sgeq 0\end{equation}
for every $G$--equivariant subsheaf $\xxF'$ of every $(G,h)$--constellation $\xxF$, both generated in $D_-$, then in consideration of Theorem \ref{equivGITstab} we would obtain that a $(G,h)$--con\-stellation is $\theta$--(semi)stable if and only if it is GIT--(semi)stable.
%By Proposition \ref{Dminusgenuegt} and the proof of the above theorem, it would even be enough to show \eqref{thequivtth} for $\xxF'$ generated in $D_-$.  

The equivalence \eqref{thequivtth} might be asking too much for, but in the following subsection we show at least that $\theta$--stability implies GIT--stability (Theorem \ref{thetaimpliesGIT}). As the theorem suggests, we therefore compare $\theta$ and $\widetilde \theta$ and we show that $\theta$--stability implies $\widetilde \theta$--stability.

\subsection{Comparison of $\theta$ and $\widetilde\theta$} \label{comparison}

Let $\xxF$ be any $(G,h)$--constellation and $\xxF' \subset \xxF$ a $G$--equivariant coherent subsheaf. We have defined two functions on $\Coh^G(X)$:
\begin{align*}
\theta(\xxF') \ &= \ \sum_{\rho \in D_-}\theta_{\rho}h'(\rho) \quad + \sum_{\sigma \in D\setminus D_-}\theta_{\sigma}h'(\sigma) \quad + \sum_{\tau \in \Irr G \setminus D}\theta_{\tau}h'(\tau),\\
%\intertext{and}
 \widetilde \theta(\xxF') \ &=\ \sum_{\rho \in D_-}\Big(\kappa_{\rho} + \chi_{\rho} -\frac{\kappa(\xxF)}{\dim A}\Big)h'(\rho) \quad + \sum_{\sigma \in D\setminus D_-}\kappa_{\sigma}h'(\sigma). 
\end{align*}
The main difference is that $\theta$ is defined as the sum over infinitely many elements while the number of summands in $\widetilde \theta$ is finite. We define the  part outside $D$ of $\theta$ by
\begin{equation} \label{def_SD}
S_D=S_D(h) := \sum_{\tau \in \Irr G \setminus D}\theta_{\tau}h(\tau).
\end{equation}

To compare $\theta$ and $\widetilde \theta$ we make the following approach for choosing the character $\chi$ and the weights $\kappa$ in the definition of our ample line bundle $\xLl$:
\begin{equation} \label{choicechikappa}
\left\{
    \begin{array}{llll}  
\kappa_{\rho} &\in \ratQ_{>0} \quad \text{ arbitrary} &&\text{for } \rho \in D \cap (D_- \cup D_0),\\
\kappa_{\rho} &= \theta_{\rho} + \frac{S_D}{d \cdot h(\rho)} &&\text{for } \rho \in D \cap D_+,\\
\chi_{\rho} &= \theta_{\rho} - \kappa_{\rho} + \frac{\kappa(\xxF)}{\dim A}  &&\text{for } \rho \in D_-,    
\end{array}
\right.
\end{equation}
where $d := \# (D\cap D_+)$ is the number of non-zero summands in the second sum in the definition of $\widetilde \theta$. We recall that if $\rho \in D_- \cup D_+$, then $h(\rho) \neq 0$ by assumption on $D_0$ (see the end of Section \ref{deff}). The condition $\theta(\xxF) = 0$ implies $S_D = - \sum_{\sigma \in D} \theta_{\sigma}h(\sigma) \in \ratQ$, so $\kappa_{\rho} \in \ratQ$.
Since $D \supset D_-$, the inequality $S_D \geq 0$ holds, and thus $\kappa_\rho>0$ for all $\rho \in D$.

Let us also note the following facts substantiating why the choice \eqref{choicechikappa} for $\chi$ and $\kappa$ is natural:
\begin{itemize}
\item Since $\theta_{\rho} < 0$ and $\kappa_{\rho} > 0$ for every $\rho \in D_-$, we automatically have
$$\chi_{\rho} = \theta_{\rho} - \kappa_{\rho} + \frac{\kappa(\xxF)}{\dim A} < \frac{\kappa(\xxF)}{\dim A},$$
so the prerequisites of Lemma \ref{GITssiso} and Theorem \ref{equivGITstab} are always satisfied.
\item One easily calculates $\sum_{\rho \in D_-}\chi_{\rho}h(\rho) = \theta(\xxF)$, and thus if $\xxF$ is $\theta$--semistable we obtain $\sum_{\rho \in D_-}\chi_{\rho}h(\rho)=0$.
\item Let $\xxF$ be a $(G,h)$--constellation. For any $G$--equivariant coherent subsheaf $\xxF'$ of $\xxF$, plugging in \eqref{choicechikappa} in Definition \ref{thetatilde} gives
\begin{equation}  \label{eqthetatilde}
\widetilde \theta(\xxF')\ =\ \sum_{\rho \in D} \theta_{\rho}h'(\rho) \ + \ \frac{S_D}{d}\sum_{\sigma \in D \cap D_+}\frac{h'(\sigma)}{h(\sigma)},
\end{equation}
in particular $\widetilde \theta(\xxF) = \theta(\xxF)$.
\end{itemize}

\begin{remark}
If the support $D_- \cup D_+$ of $\theta$ is finite, then one may take $D \supset D_- \cup D_+$. In this case the summand $S_D$ vanishes and \eqref{eqthetatilde} yields
$$\widetilde \theta(\xxF') = \sum_{\rho \in D_- \cup D_+} \theta_{\rho}h'(\rho) = \theta(\xxF'), \text{ for all } \xxF' \subset \xxF.$$
%In particular, if $G$ is a finite group, $\theta$--(semi)stability and GIT--(semi)stability coincide as in the construction of Craw and Ishii \cite{CI:2004}.
%But for a reductive group $G$, the support of $\theta$ will be infinite in general for otherwise the $(G,h)$--constellations which are $\theta$--semistable but not $\theta$--stable might not be generated in $D_-$.
\end{remark}
%\goodbreak

For comparing $\theta$ to $\widetilde \theta$, we consider $\widetilde \theta = \widetilde \theta_D$ when the finite subset $D \subset \Irr G$ varies. We obtain the following error terms:

\begin{prop}\label{errorestimate}
If $\widetilde D \supset D $, then for any $G$--equivariant coherent subsheaf $\xxF'$ of a $(G,h)$--con\-stellation $\xxF$ we have
\begin{align*}
| \widetilde \theta_{\widetilde D}(\xxF') - \widetilde \theta_{D}(\xxF')| & \leq \sum_{\tau \in \widetilde D \setminus D } \left( \theta_{\tau}h(\tau) + \frac{S_{\widetilde D}}{\widetilde d} \right),\\
\intertext{where $S_{\widetilde D}$ is defined by \eqref{def_SD} and $\widetilde d := \# (\widetilde D \cap D_+)$. Further, we have}
| \theta(\xxF') - \widetilde \theta_{D}(\xxF')|  & \leq \sum_{\tau \in \Irr G \setminus D} \theta_{\tau}h(\tau).
\end{align*}
\end{prop}

\begin{proof}
Using \eqref{eqthetatilde}, we write 
\begin{align*}
\widetilde \theta_D(\xxF') &= \sum_{\rho \in D_-} \theta_{\rho}h'(\rho) \;\; + \sum_{\sigma \in D \cap D_+}\theta_{\sigma}h'(\sigma) \; + \; \frac{S_D}{d}\sum_{\sigma \in D \cap D_+}\frac{h'(\sigma)}{h(\sigma)} \qquad \text{and}\\
\widetilde \theta_{\widetilde D}(\xxF') &= \sum_{\rho \in D_-} \theta_{\rho}h'(\rho) \;\; + \sum_{\sigma \in \widetilde D \cap D_+}\theta_{\sigma}h'(\sigma) \; + \; \frac{S_{\widetilde D}}{\widetilde d}\sum_{\sigma \in \widetilde D \cap D_+}\frac{h'(\sigma)}{h(\sigma)}.
\end{align*}
In \cite[Proposition 4.3.3]{Bec:thesis}, the determination of their difference is carried out:
\begin{align*}
\widetilde \theta_{\widetilde D}(\xxF') - \widetilde \theta_{D}(\xxF') &= \sum_{\tau \in (\widetilde D \setminus D) \cap D_+} \left( \theta_{\tau}h(\tau) + \frac{S_{\widetilde D}}{\widetilde d} \right) \left( \frac{h'(\tau)}{h(\tau)} - \frac{1}{d}\sum_{\sigma \in D \cap D_+} \frac{h'(\sigma)}{h(\sigma)}\right)\\
              &-\frac{1}{d}\left(\sum_{(\widetilde D \setminus D) \cap D_0} \frac{S_{\widetilde D}}{\widetilde d} \right) \left(\sum_{ \sigma \in D \cap D_+} \frac{h'(\sigma)}{h(\sigma)} \right).
\end{align*}

Further, for every $\tau \in D_- \cup D_+$ we have $0 \leq \frac{h'(\tau)}{h(\tau)} \leq 1$, so it follows \\ $\Big|\frac{h'(\tau)}{h(\tau)} - \frac{1}{d}\sum_{\sigma \in D \cap D_+} \frac{h'(\sigma)}{h(\sigma)}\Big| \leq 1$. We deduce
\begin{align*}
| \widetilde \theta_{\widetilde D}(\xxF') - \widetilde \theta_{D}(\xxF') | & \leq \sum_{\tau \in (\widetilde D \setminus D) \cap D_+} \left| \theta_{\tau}h(\tau) + \frac{S_{\widetilde D}}{\widetilde d} \right|+\sum_{\tau \in (\widetilde D \setminus D) \cap D_0} \left| \frac{S_{\widetilde D}}{\widetilde d} \right|\\
&=\sum_{\tau \in (\widetilde D \setminus D) } \left( \theta_{\tau}h(\tau)+ \frac{S_{\widetilde D}}{\widetilde d} \right).
\end{align*}

A similar computation gives the second upper bound.
\end{proof}

The set $\xxD = \{ D \subset \Irr G \mid D \supset D_-\}$ is directed with respect to inclusion. In this sense, we can take the limit over these sets.
This allows us to reveal the relation between $\theta$ and $\widetilde \theta$:

\begin{cor}\label{limit}
The function $\theta$ is the pointwise limit of the functions $\widetilde \theta_{D}$ as $D$ converges to $\Irr G$:
$$\theta(\xxF') = \lim_{D \in \xxD} \widetilde \theta_D (\xxF'), \text{ for all }\xxF' \subset \xxF.$$
\end{cor}

\begin{proof}
Since $\theta(\xxF) = \sum_{\tau \in \Irr G} \theta_{\tau}h(\tau)$ is convergent, the sum $\sum_{\tau \in \Irr G \setminus D} \theta_{\tau}h(\tau)$ converges to $0$ when $D$ becomes larger. Then the result follows from the second inequality of Proposition \ref{errorestimate}.
\end{proof}

In general, equality will only hold in the limit, but not for finite $D$. We use this corollary in the next proposition to show that every $\theta$--stable $(G,h)$--constellation is also $\widetilde \theta$--stable.

\begin{prop}\label{chooseD}
There is a finite subset $D \subset \Irr G$ such that the following holds: If $\xxF$ is a $\theta$--stable $(G,h)$--constellation and $\xxF'$ a non-zero proper $G$--equivariant subsheaf of $\xxF$ generated in $D_-$, then for every finite set $\widetilde D$ containing $D$ we have $\widetilde \theta_{\widetilde D}(\xxF') > 0$.
\end{prop}

\begin{proof}
By Proposition \ref{endlvieleHilbertfkt}, the set 
$$\{\theta(\xxF'')\mid \xxF'' \subset \xxF \text{ a non-zero proper $G$--equivariant subsheaf generated in } D_-\}$$
is finite. Let $\theta_0>0$ be its minimum. 
%In particular, we have $\theta(\xxF') \geq \theta_0$ for every non-zero proper $G$--equivariant subsheaf $\xxF' \subset \xxF$ by Proposition \ref{Dminusgenuegt}.
If we fix $\varepsilon > 0$, by Corollary \ref{limit} there is a subset $D = D(\varepsilon) \subset \Irr G$ such that $| \theta(\xxF') - \widetilde \theta_{\widetilde D}(\xxF') | < \varepsilon$ for every $\widetilde D \supset D$. 
Now if we choose $\varepsilon < \theta_0$, we obtain $D = D(\varepsilon)$ such that for every $\widetilde D \supset D$ we have 
$$\widetilde \theta_{\widetilde D}(\xxF') >  \theta(\xxF') - \varepsilon  \geq \theta_0 - \varepsilon > 0.$$
\end{proof}

Now we summarize:
\begin{thm}\label{thetaimpliesGIT}
Let $\theta \in \ratQ^{\Irr G}$ be a stability condition on the set of $(G,h)$--con\-stellations on $X$ with $D_-$ finite and $\langle\theta,h\rangle = 0$, and
let $\kappa$ and $\chi$ be defined as in \eqref{choicechikappa}. For $\xxH := \bigoplus_{\rho \in D_-} \kxC^{h(\rho)} \tensor V_{\rho} \tensor \ssh_X$ let us consider the invariant Quot scheme $\xxGQuotH$ and the ample line bundle $\xLl = \bigotimes_{\sigma \in D} (\det \mathcal{U}_{\sigma})^{\kappa_{\sigma}}$ on $\xxGQuotH$ with $D \subset \Irr G$ large enough in the sense of Proposition \ref{chooseD}. 
Let the natural $\Gamma$--linearization on $\xLl$ be twisted by the character $\chi$. 
For $\xxF' \in \Coh^G(X)$ let us set
$\widetilde \theta(\xxF') = \sum_{\rho \in D_-}\big(\kappa_{\rho} + \chi_{\rho} -\frac{\kappa(\xxF)}{\dim A}\big)h'(\rho) + \sum_{\sigma \in D\setminus D_-}\kappa_{\sigma}h'(\sigma)$. With these choices of $D$, $\kappa$, $\chi$ and $\widetilde \theta$, every $\theta$--stable $(G,h)$--constellation is $\widetilde \theta$--stable and hence GIT--stable.
\end{thm}
\begin{proof}
This is a direct consequence of Proposition \ref{chooseD} and Theorem \ref{equivGITstab}.
\end{proof}

On the level of functors, we obtain the following:

\begin{cor}
With the same notation and choices as in Theorem \ref{thetaimpliesGIT}, the moduli functor $\fsMth$ of $\theta$--stable $(G,h)$--constellations on $X$ is a subfunctor of the moduli functor $\fsMck$ of GIT--stable $(G,h)$--constellations on $X$.
\end{cor}

\begin{remark}
Even if we will not use this fact in the sequel, it is worth noticing that $\widetilde \theta$--semistability implies $\theta$--semistability. Indeed, if there exists $D$ such that for every $D \subset D'$ we have $\xxF$ is $\widetilde \theta_{D'}$--semistable, then it follows from Corollary \ref{limit} that $\xxF$ is also $\theta$--semistable. Therefore the moduli functor $\fMck$ of GIT--semistable $(G,h)$--constellations on $X$ is a subfunctor of the moduli functor $\fMth$ of $\theta$--semistable $(G,h)$--constellations on $X$.
\end{remark}

% *** 5 ***************************************************************************************************************************************
\section{The moduli space of $\theta$--stable $(G,h)$--constellations}\label{modspace}

In this section, we use the notation and assumptions of Theorem \ref{thetaimpliesGIT}. The preceding sections leave us with the following situation: We have
\begin{align*}
&\xxGQuotH^{ss} &&\hspace{-1em}:= \left\{[q\colon \xxH \twoheadrightarrow \xxF] \in \xxGQuotH \middle| [q] \text{ is GIT--semistable}\right\},\\
&\xxGQuotH^s    &&\hspace{-1em}:= \left\{[q\colon \xxH \twoheadrightarrow \xxF] \in \xxGQuotH \middle| [q] \text{ is GIT--stable}\right\},\\
&\xxGQuotH^{ss}_{\theta} &&\hspace{-1em}:=\left\{[q\colon \xxH \twoheadrightarrow \xxF] \in \xxGQuotH \middle| \begin{array}{l}\xxF \text{ is } \theta\text{--semistable and } \\ $[q]$ \text{ originates from $\xxF$}\end{array}\right\},\\
&\xxGQuotH^s_{\theta} &&\hspace{-1em}:= \left\{[q\colon \xxH \twoheadrightarrow \xxF] \in \xxGQuotH \middle| \begin{array}{l} \xxF \text{ is } \theta\text{--stable and }\\ $[q]$ \text{ originates from $\xxF$}\end{array}\right\}
\end{align*}
and set-theoretic inclusions
$$\xxGQuotH^s_{\theta} \subset \xxGQuotH^s  \subset \xxGQuotH^{ss} \subset \xxGQuotH^{ss}_{\theta}.$$
Forgetting the choice of the particular quotient map, this yields inclusions
\begin{small}
\begin{align*}
\left\{\begin{array}{c}\theta\text{--stable }\\(G,h)\text{--constellations}\end{array}\right\} 
&\subset \left\{\begin{array}{c}\text{GIT--stable }\\(G,h)\text{--constellations}\end{array}\right\} \\
&\subset \left\{\begin{array}{c}\text{GIT--semistable }\\(G,h)\text{--constellations}\end{array}\right\}
\subset \left\{\begin{array}{c} \theta\text{--semistable }\\(G,h)\text{--constellations}\end{array}\right\}.
\end{align*}
\end{small}
On the level of functors this translates into a sequence 
$$\fsMth \subset \fsMck \subset \fMck \subset \fMth.$$

On the level of schemes we show that $\xxGQuotH^s_{\theta}$ is an open subscheme of $\xxGQuotH^s$ in Section \ref{openness}

In Section \ref{corepr} we prove that $\fMck$, $\fsMck$ and $\fsMth$ are corepresented by the categorical quotient $\xxGQuotH^{ss} \red_{\xLl_{\chi}} \Gamma$, the geometric quotient $\xxGQuotH^s /\Gamma$ and the geometric quotient $\sMth := \xxGQuotH^s_{\theta} /\Gamma$, respectively. Thus we obtain
$$ \begin{xy} \xymatrix{
\hspace{-1em}\xxGQuotH^s_{\theta}\hspace{-1em} \ar@{->>}[d] & \hspace{-3em} \subset \hspace{-3em} & \hspace{-1em}\xxGQuotH^s\hspace{-1em} \ar@{->>}[d] & \hspace{-3em} \subset \hspace{-3em} & \hspace{-1em}\xxGQuotH^{ss}\hspace{-1em} \ar@{->>}[d]\\
\hspace{-1em}\sMth\hspace{-1em} & \hspace{-3em} \subset \hspace{-3em} & \hspace{-1em}\xxGQuotH^s/\Gamma\hspace{-1em} & \hspace{-3em} \subset \hspace{-3em} & \hspace{-1em}\xxGQuotH^{ss}\red_{\xLl_{\chi}}\Gamma.\hspace{-1em}
}\end{xy}$$
We deduce that $\sMth$, the moduli space of $\theta$--stable $(G,h)$--constellations, is an open subscheme of $\xxGQuotH^s /\Gamma$ and is therefore quasiprojective. It generalizes the invariant Hilbert scheme as we have shown in Section \ref{HilbgleichMth}.

In Section \ref{quotmap}, under the extra assumption that $h(\rho_0)=1$, we conclude the construction of $\sMth$ as a moduli space over the quotient $X \red G$: we construct a morphism $\sMth \to X\red G$ generalizing the Hilbert--Chow morphism.

\subsection{Openness of $\theta$--stability}\label{openness}

In this section we prove that $\xxGQuotH^s_{\theta}$ is an open subscheme of $\xxGQuotH^s$.

\begin{prop}\label{openinfamilies}
Being $\theta$--(semi)stable is an open property in flat families of $(G,h)$--constellations.
\end{prop}
\begin{proof}
We proceed analogously to \cite[Proposition 2.3.1]{HL:2010}. Let $f \colon \mathcal{X} \to S$ be a family of affine $G$--schemes and $\xxFf$ a flat family of $(G,h)$--constellations on $\mathcal{X}$. Let
\begin{align*}
&H &&\hspace{-1.2em}:= \left\{h'' \text{ a Hilbert function } \middle| \begin{array}{l}\exists \; s \in S \text{ and a surjection } \alpha(s)\colon \xxFf(s) \to \xxF''\\ \text{with } \ker \alpha(s) \text{ generated in } D_- \text{ and } h_{\xxF''} = h''\end{array}\right\}, \\
&H^{ss} &&\hspace{-1.2em}:= \{h'' \in H \mid \langle\theta,h''\rangle > 0\}, \\
&H^s &&\hspace{-1.2em}:= \{h'' \in H \mid \langle\theta,h''\rangle \geq 0\}.
\end{align*}
By Proposition \ref{endlvieleHilbertfkt} and Remark \ref{thetaquot}, $H$ is finite. For each Hilbert function $h''$ in $H$ we consider the relative invariant Quot scheme $\pi_{h''}\colon \relGQuot{\mathcal{X}/S}{\xxFf}{h''} \to S$ with fibres $\xxGQuot{\xxFf(s)}{h''}$ over $s \in S$. Since the multiplicities of the $\xxFf(s)$ are finite, the map $\pi_{h''}$ is projective by \cite[Proposition B.5]{Bec:thesis}. Thus its image is a closed subset of $S$.
Remark \ref{thetaquot} says that $\xxFf(s)$ is $\theta$--(semi)stable if and only if the Hilbert function $h''$ of every strict quotient of $\xxFf(s)$ satisfies $\langle\theta,h''\rangle \sleq 0$. Accordingly, $\xxFf(s)$ is $\theta$--(semi)stable if and only if $s$ is not contained in the finite, hence closed, union $\bigcup_{h'' \in H^{(s)s}} \imag(\pi_{h''})$.
\end{proof}

\begin{cor}  \label{open_sub}
$\xxGQuotH^s_{\theta}$ is an open subscheme of $\xxGQuotH^s$.
\end{cor}

\begin{proof}
Let 
$$\xxFf \in \fGQuotH \big(\xxGQuotH \big)$$ 
be the universal family over $\xxGQuotH$, so that the fibre $\xxFf(\xxF)$ equals $\xxF$. Since by Proposition \ref{openinfamilies} the property of being $\theta$--stable is open in flat families, the set
$$\Xi:=\{[q\colon \xxH \twoheadrightarrow \xxF] \in \xxGQuotH \mid \xxFf(\xxF) \text{ is } \theta\text{--stable}\}$$
is open in $\xxGQuotH$. 
On the other hand, by Theorem \ref{thetaimpliesGIT} we know that  $\xxGQuotH_{\theta}^s = \Xi \cap \xxGQuotH^s$, whence the result.
\end{proof}

Since $\xxGQuotH^s_{\theta}$ is an open subscheme of $\xxGQuotH$, it is a quasiprojective scheme.

\subsection{Corepresentability}\label{corepr}

Let $R \!:= \xxGQuotH^{ss}$, $R^s \!:= \xxGQuotH^{s}$ and $R^s_{\theta} \!:= \xxGQuotH^{s}_{\theta}$ be the open subschemes of $\xxGQuotH$ corresponding to  GIT--semistable, GIT--stable and $\theta$--stable quotients of $\xxH$, respectively. 
For elements $[q\colon \xxH \twoheadrightarrow \xxF] \in \xxGQuotH$, the sheaf $\xxF = q(\xxH) = q(\bigoplus_{\rho \in D_-} A_{\rho}\tensor V_{\rho} \otimes_{\kxC} \ssh_X)$ is automatically generated in $D_-$. 
Moreover, for $[q\colon \xxH \twoheadrightarrow \xxF] \in R$ the maps $\varphi_{\rho}\colon A_{\rho} \to \xxF_{\rho}$, $a \mapsto (v \mapsto q(a \otimes v \otimes 1))$ are isomorphisms for every $\rho \in D_-$ by Lemma \ref{GITssiso} and since the inequality $\chi_{\rho} < \frac{\kappa(\xxF)}{\dim A}$ holds.
As presented in Subsection \ref{parameters}, the choice of these isomorphisms is described by the action of the group $\Gamma$ defined by \eqref{def Gamma}, which acts (with finite stabilizers) on $\xxGQuotH$ from the right by left multiplication on the components of $\xxH$. The subsets $R$ and $R^s$ are invariant under this action, and the same holds for $R^s_{\theta}$ since the action of an element in $\Gamma$ does not change $\xxF$. %the Hilbert function of the subsheaves of $\xxF$.

Here comes the main result of this paper:
%To deal with the ambiguity of the choice of the $\varphi_{\rho}$, we have the following relation between the moduli problem and quotients of this group action:

\begin{thm}
With the notation above, the functors $\fMck$, $\fsMck$ and $\fsMth$ of Definition \ref{thetafunctors} and Definition \ref{GITfunctors} are corepresented by the categorical quotient $R\red_{\xLl_{\chi}}\Gamma$, the geometric quotient $R^s/\Gamma$ and the geometric quotient $R^s_{\theta}/\Gamma$, respectively.  
\end{thm}

\begin{proof} We proceed analogously to \cite[Lemma 4.3.1]{HL:2010}.\\\hspace{-2.5em}
The quotients by $\Gamma$ and $\Gamma'=\prod_{\rho \in D_-} Gl(A_\rho)$ coincide since multiples of the identity act trivially. Hence we can consider the action of $\Gamma'$ instead of the one of $\Gamma$ in the sequel.
Let $S$ be a noetherian scheme over $\kxC$ and $\xxFf$ a flat family of $(G,h)$--con\-stellations generated in $D_-$ which is parametrized by $S$, so that for every $s \in S$ the fibre $\xxFf(s)$ is a $(G,h)$--constellation on X. Let $p\colon X \times S \to S$ denote the projection.  We look at the isotypic decomposition 
$$p_*\xxFf \cong \bigoplus_{\rho \in \Irr G} \xxFf_{\rho}\tensor V_{\rho}.$$
The conditions that $G$ is reductive, $p$ is affine and $\xxFf$ is flat over $S$ yield that the $\xxFf_{\rho}$ are locally free $\ssh_S$--modules of rank $h(\rho)$ and that we have $(\xxFf_{\rho})(s) = \xxFf(s)_{\rho}$.
We define the $\ssh_S$--submodule
\begin{equation}\label{VF}
V_{\xxFf}^- := \bigoplus_{\rho \in D_-}(p_*\xxFf)_{(\rho)} = \bigoplus_{\rho \in D_-} \xxFf_{\rho} \tensor V_{\rho} \subset p_*\xxFf.
\end{equation}
The pullback of the inclusion $i\colon V_{\xxFf}^- \hookrightarrow p_*\xxFf$ composed with the natural surjection $\alpha\colon p^*p_*\xxFf \twoheadrightarrow \xxFf$ corresponding to the identity under the adjunction \linebreak $\xxHom(p^*p_*\xxFf,\xxFf)$ $\cong \xxHom(p_*\xxFf,p_*\xxFf)$ yields a morphism
$$\varphi_{\xxFf} := \alpha \circ p^*i\colon p^*V_{\xxFf}^- \to \xxFf.$$
Fiberwise, $\varphi_{\xxFf}(s)\colon (p^*V_{\xxFf}^-)(s) \tensor \ssh_X \to \xxFf(s)$ % (p^*V_{\xxFf}^-)(s) = \bigoplus_{\rho \in D_-} (\xxFf(s))_{\rho} \tensor V_{\rho}
is surjective since each $\xxFf(s)$ is generated in $D_-$ as an $\ssh_X$--module. So $\varphi_{\xxFf}$ is also surjective.

Let $A_V := \bigoplus_{\rho \in D_-} A_{\rho} \tensor V_{\rho}$ and $\pi\colon \mathbbm{I}(\xxFf) := \Isom_G(A_V \tensor \ssh_{S}, V_{\xxFf}^-) \to S$ the $G$--equivariant frame bundle associated to $V_{\xxFf}^-$ as described in \cite[Appendix A]{Bec:thesis}. It parametrizes $G$--equivariant isomorphisms $A_V \tensor \ssh_{S} \to V_{\xxFf}^-$ and gives us a canonical morphism $\alpha\colon A_V \tensor \ssh_{\mathbbm{I}(\xxFf)} \to \pi^*V_{\xxFf}^-$. 

Now we consider the product $\pi_X := \xid_X \times \pi \colon X \times \mathbbm{I}(\xxFf) \to X \times S$ and the universal trivialization $\alpha \tensor \xid_X \colon A_V \tensor \ssh_{X\times \mathbbm{I}(\xxFf)} = \xxH \tensor \ssh_{\mathbbm{I}(\xxFf)} \to (p \circ \pi_X)^*V_{\xxFf}^-$ 
on $X \times \mathbbm{I}(\xxFf)$. Thus we obtain a canonically defined quotient
$$[\pi_X^*\varphi_{\xxFf}\circ (\xid_X \tensor \alpha)\colon \xxH \tensor \ssh_{\mathbbm{I}(\xxFf)} \to \pi_X^*p^*V_{\xxFf}^- \to \pi_X^*\xxFf] \in \fGQuotH(\mathbbm{I}(\xxFf)),$$
which in turn yields a classifying morphism
$$\phi_{\xxFf}\colon \mathbbm{I}(\xxFf) \longrightarrow \xxGQuotH,\; \psi \longmapsto [q_{\psi}\colon \xxH \twoheadrightarrow (\pi_X^*\xxFf)(\psi) = \xxFf(\pi(\psi))].$$ 
As for the Quot scheme, the gauge group $\Gamma'$ acts on $\mathbbm{I}(\xxFf)$ from the right. Here, $\pi\colon \mathbbm{I}(\xxFf) \to S$ is even a principal $\Gamma'$--bundle. 
By construction, $\phi_{\xxFf}$ is $\Gamma'$--equivariant and we have $\phi_{\xxFf}^{-1}(R) = \pi^{-1}(S^{ss})$, where $S^{ss} = \{s \in S \mid \xxFf(s) \text{ GIT--semistable}\}$. If $S$ parametrizes GIT--semistable sheaves, we even have $\phi_{\xxFf}^{-1}(R) = \pi^{-1}(S) = \mathbbm{I}(\xxFf)$, hence $\phi_{\xxFf}(\mathbbm{I}(\xxFf)) = \phi_{\xxFf}(\phi_{\xxFf}^{-1}(R)) \subset R$. This means that in fact we have $\phi_{\xxFf}\colon \mathbbm{I}(\xxFf) \to R$. This morphism induces a transformation of functors
$$\underline{\mathbbm{I}(\xxFf)}/\underline{\Gamma'} \to \underline{R}/\underline{\Gamma'}.$$
Since $\pi\colon \mathbbm{I}(\xxFf) \to S$ is a principal $\Gamma'$--bundle, $S$ is a categorical quotient of $\mathbbm{I}(\xxFf)$, 
so that we obtain an element in $(\underline{R}/\underline{\Gamma'})(S)$. 
Thus we have constructed a transformation $\fMck \to \underline{R}/\underline{\Gamma'}$. 

Denoting $p_R\colon X\times R \to R$, the universal family $[q\colon p_R^*\xxH \twoheadrightarrow \mathcal{U}]$ on $R$ yields an inverse by mapping the classifying morphism $(\xi\colon S \to R) \in(\underline{R}/\underline{\Gamma'})(S)$ to the family $(\xid_X \times \xi)^*\mathcal{U} \in \fMck(S)$. % wie genau?

Altogether this means that a scheme $M$ corepresents $\fMck$ if and only if it corepresents $\underline{R}/\underline{\Gamma'}$, hence if and only if it is a categorical quotient of $R$ by $\Gamma'$.

The same proof literally goes through replacing GIT--semistability by GIT--sta\-bility and $R$ and $\fMck$ by $R^s$ and $\fsMck$, respectively, as well as replacing GIT--semistability by $\theta$--stability and $R$ and $\fMck$ by $R^s_{\theta}$ and $\fsMth$, respectively. Finally, $R^s/\Gamma'=R^s/\Gamma$ is even a geometric quotient, and $R_\theta^s$ is a $\Gamma'$-invariant (open) subset of $R^s$, hence $R_\theta^s/\Gamma'=R_\theta^s/\Gamma$ is also a geometric quotient. 
\end{proof}

\begin{definition}
The scheme $\sMth := \xxGQuotH^s_{\theta}/\Gamma$ is called the \textit{moduli space of $\theta$--stable $(G,h)$--constellations}.
We denote by $\Mth$ the closure of $\sMth$ in $\xxGQuotH^{ss}\red_{\xLl_{\chi}}\Gamma$.
\end{definition}

By \cite[\S 0.2 (4)]{GIT:1994}, the quotient map $\nu\colon \xxGQuotH^s \to \xxGQuotH^s/\Gamma$ is open. Thus, by Corollary \ref{open_sub}, the image $\sMth = \nu(\xxGQuotH^s_{\theta})$ is open in $\xxGQuotH^s/\Gamma$; in particular $\sMth$ and $\Mth$ are quasiprojective schemes. 

\begin{remark}
In Section \ref{HilbgleichMth} we have already seen that if $h(\rho_0) = 1$ and
if $\theta$ is chosen such that $D_- = \{\rho_0\}$ we recover the invariant Hilbert scheme: 
$$\sMth = \invHilb{G}{h}{X}.$$
\end{remark}

\subsection{The map to the quotient $X \red G$} \label{quotmap}

As before, let $\xxH=\bigoplus_{\rho \in D_-} A_\rho \otimes_\kxC V_\rho \otimes_\kxC \ssh_X$ and $h: \Irr G \to \natN$ a Hilbert function. 
For the invariant Quot scheme, Jansou \cite[\S 1.4]{Jan:2006} constructed an analogue of the Hilbert--Chow morphism
$$\gamma\colon \xxGQuotH \longrightarrow \Quot(\xxH^G,h(\rho_0)),\;
[q\colon \xxH \twoheadrightarrow \xxF] \longmapsto [q|_{\xxH^G}\colon \xxH^G \twoheadrightarrow \xxF^G].$$
In the case where $h(\rho_0)=1$, we show how $\gamma$ induces a morphism $\Mth \to X/\!/G$. 

%In the case where $h(\rho_0) = 1$, we extend the restriction $\gamma|_{\xxGQuotH^{ss}}$ to a morphism to $X \red G$:
\begin{thm}\label{genHC}
If $h(\rho_0) = 1$, there is a morphism $\xxGQuotH^{ss} \to X \red G$, which yields a morphism 
$$\eta\colon \Mth \to X\red G,\; \xxF \mapsto \supp \xxF^G.$$
\end{thm}

\begin{proof}
Let $S$ be a noetherian scheme over $\kxC$ and $[q\colon\pi^*\xxH \twoheadrightarrow \xxFf]$ an element of $\fGQuotH^{ss}(S)$, where $\pi\colon X \times S \to X$. Then we have $\gamma_S(q)\colon \ssh_S \tensor \xxH^G \to \xxFf^G$. Since every fibre $q(s)$ is GIT--semistable, the morphism $\varphi_{\rho_0}\colon A_{\rho_0} \to \xxFf^G(s)$
defined in \eqref{iso} is an isomorphism for every $s \in S$. Hence $\gamma_S(q)$ restricted to the subset $\ssh_S \tensor \ssh_X^G \cong \ssh_S \tensor A_{\rho_0} \tensor \ssh_X^G$ of $\ssh_S \tensor \xxH^G$ maps surjectively to $\xxFf^G$.
Consider the composite morphism
$$\psi\colon\ssh_S \xrightarrow{\xid \otimes 1} \ssh_S \tensor \ssh_X^G \twoheadrightarrow \xxFf^G.$$
If $\psi(s)\!:\!\ssh_S (s) \to \xxFf^G(s)$ were $0$ for some $s\in S$, the map $\ssh_S \tensor \ssh_X^G \to \xxFf^G$ would not be surjective on the fibre $\xxFf^G(s)$, so this cannot happen. Thus $\psi$ is nowhere $0$. The $\ssh_S$--modules $\ssh_S$ and $\xxFf^G$ are both locally free of rank $1$, so $\psi$ is an isomorphism. This shows that $\xxFf^G$ corresponds to a subscheme $Z \subset S \times X \red G$.
With the notation
$$ \begin{xy} \xymatrix{
Z\; \ar@{^(->}^/-3mm/{i}[r] \ar@{->}_{p}^/1mm/{\cong}[dr] & S \times X \red G \ar@{->}^/1mm/{pr_2}[r] \ar@{->}[d] & X \red G\\
 & S & 
}\end{xy}$$
we obtain a morphism
$$pr_2 \circ i \circ p^{-1}\colon S \to X \red G.$$
This procedure is compatible with base change, 
%Indeed, let $g\colon T \to S$ be a morphism of noetherian schemes over $\kxC$. This yields $[(id_X \times g)^*q\colon \pi_T^*\xxH \twoheadrightarrow (id_X \times g)^*\xxFf] \in \fGQuotH^{ss}(T)$, where $\pi_T\colon X \times T \to X$ is the projection to $X$. The invariants satisfy
%$$((id_X \times g)^*\xxFf)^G = g^*\xxFf^G \cong g^*\ssh_S \cong \ssh_T.$$
%Hence, the subscheme corresponding to $((id_X \times g)^*\xxFf)^G$ is $Z_T = Z \times_S T$ and we have
%$$ \begin{xy} \xymatrix{
%Z_T\; \ar@{^(->}^/-3mm/{j}[r] \ar@{->}_{p_T}^/1mm/{\cong}[dr] & T \times X \red G \ar@{->}[r] \ar@/_5mm/@{->}_{pr_2\circ (g \times \xid_{X \red G})}[rr] \ar@{->}[d] & S \times X \red G \ar@{->}[r] & X \red G\\
% & T & 
%}\end{xy}$$
%The above construction yields a morphism
%$$(pr_2 \circ (g \times \xid_{X \red G})) \circ j \circ p_T^{-1}\colon T \to X \red G.$$
%We have the following commuting diagram:
%$$ \begin{xy} \xymatrix{
%T \times X \red G \ar@{->}^{g \times \xid_{X \red G}}[d] & \;Z_T \ar@{->}[d] \ar@{_(->}_/-3mm/{j}[l] \ar@{->}^{p_T}[r] & T \ar@{->}^{g}[d]\\
%S \times X \red G & \;Z \ar@{_(->}_/-3mm/{i}[l] \ar@{->}^{p}[r] & S 
%}\end{xy}$$
%In particular, we have $i \circ p^{-1} \circ g = (g \times \xid_{X \red G}) \circ j \circ p_T^{-1}$, so that the morphism for $T$ is the composition of the morphism for $S$ with $g$:
%$$ \begin{xy} \xymatrix@M=0.7em{
%T \ar@{->}^{g}[d] \ar@{->}^/-2mm/{(pr_2 \circ (g \times \xid_{X \red G})) \circ j \circ p_T^{-1}}[rrr] & & & X \red G\\
%S \ar@{->}_{pr_2 \circ i \circ p^{-1}}[urrr] 
%}\end{xy}$$
so we have constructed a morphism of functors
$$\fGQuotH^{ss} \to \Mor(\cdot, X\red G).$$
Plugging in $\xxGQuotH^{ss}$, this gives a morphism of schemes
$$\eta\colon \xxGQuotH^{ss} \to X\red G.$$
By construction, for a point $[q\colon \xxH \twoheadrightarrow \xxF] \in \xxGQuotH^{ss}$ the subscheme of $X\red G$ corresponding to $\xxF^G = \ssh_X^G/I_{\xxF}$ is just its support 
$$\supp \xxF^G = \{\mathfrak{p} \in \ssh_X^G \mid \mathfrak{p} \supset I_{\xxF}\} = \left\{\sqrt{I_{\xxF}}\right\}.$$
It only consists of one point since $\dim \xxF^G = h(\rho_0) = 1$.

\noindent As $\xxF^G$ does not depend on the choice of a basis of $\xxH$, the morphism $\eta$ is $\Gamma$--in\-variant. Hence it descends to $\xxGQuotH^{ss}\red_{\xLl_{\chi}}\Gamma$. Restricting it to $\Mth$ we eventually obtain a morphism
$$\eta\colon \Mth \to X\red G,\; \xxF \mapsto \supp \xxF^G.$$
\end{proof}

\begin{remark}
Under the conditions of the theorem, the restricted morphism \linebreak $\xxGQuotH^{s} \to X \red G$ descends to $\xxGQuotH^s/\Gamma$ and its restriction to $\sMth$ yields $\eta\colon \sMth \to X\red G,\; \xxF \mapsto \supp \xxF^G$ in the same way.
\end{remark}

Thus when $h(\rho_0) = 1$ we have constructed an analogue of the Hilbert--Chow morphism for $\sMth$ and $\Mth$, which relates these moduli spaces to the quotient $X \red G$.

%\begin{remark}
%In Proposition \ref{pQuotimm} we constructed morphisms
%$$\gamma_{\rho} \colon \xxGQuotH \longrightarrow \Quot(\xxH_{\rho},h(\rho)), \; [q\colon \xxH\to \xxF] \longmapsto [q|_{\xxH_{\rho}}\colon \xxH_{\rho} \to \xxF_{\rho}],$$
%where $\gamma_{\rho_0}$ is the Hilbert--Chow morphism $\gamma$. Therefore one may adopt the proof of Theorem \ref{genHC} to this more general situation and obtain morphisms
%$$\eta_{\rho}\colon \Mth \longrightarrow \sym^{h(\rho)}(X\red G), \; \xxF \longmapsto \supp \xxF_{\rho}$$
%for an arbitrary Hilbert function $h$ and for every $\rho \in \Irr G$.
%\end{remark}

% *** 6 ***************************************************************************************************************************************
\section{Outlook}  \label{secO}

In this paper we have constructed the moduli space $\sMth$ of $\theta$--stable $(G,h)$--constellations and, when $h(\rho_0) = 1$ for $\rho_0$ the trivial representation, a morphism $\eta\colon \sMth \to X \red G$. Examples of these moduli spaces are given by invariant Hilbert schemes. The determination of further examples would be interesting in order to get an idea of the properties of these moduli spaces, e.g.~concerning smoothness, connectedness and, for symplectic varieties $X\red G$, symplecticity of $\sMth$. 
Moreover, some questions concerning the closure of $\sMth$ and the properties of $\eta$ still have to be investigated. 

Here we discuss some ideas which are worth pursuing in the future.

\subsection{The geometric meaning of points in $\Mth$}
We defined the moduli space $\Mth$ as the closure of $\sMth$ in the categorical quotient $\xxGQuotH^{ss}\red_{\xLl_{\chi}}\Gamma$ without explicitly describing its elements geometrically. A natural question is

\begin{ques}
Does the scheme $\Mth$ corepresent the moduli functor $\fMth$ of $\theta$--semistable $(G,h)$--constellations?
\end{ques}

First of all, one has to face the question if every $\theta$--semistable $(G,h)$--constellation is also GIT--semistable. Secondly, it would be interesting to determine the values of $\theta$ for which the notions of $\theta$--stability and $\theta$--semistability coincide. In this case we obtain $\Mth = \sMth$.
For example this is true for the invariant Hilbert scheme: Since $h(\rho_0) = 1$, a $(G,h)$--constellation has no non-zero proper subsheaf generated in $D_-$. Hence every $\theta$--semistable $(G,h)$--constellation is $\theta$--stable.

%any subsheaf $\xxF'$ of a $(G,h)$--constellation $\xxF$ has a Hilbert function $h'$ with $h'(\rho_0) = 0$ or $h'(\rho_0) = 1$. In the first case, $\theta(\xxF')$ is strictly positive and in the second case, $\xxF' = \xxF$ by Section \ref{HilbgleichMth}. 

In the construction of Craw and Ishii \cite{CI:2004} and King \cite{Kin:1994} the stability condition $\theta$ only consists of finitely many components. In their case, $\theta$--semistability and GIT--semistability are even equivalent, as well as $\theta$--stability and GIT--stability. It would be interesting to know if this also holds in our case. With regard to Theorem \ref{equivGITstab} this question is equivalent to the following one:

\begin{ques} \label{questionequiv}
Let $\xxF$ be a $(G,h)$--constellation and $\xxF'$ a $G$--equivariant coherent subsheaf of $\xxF$, both generated in $D_-$. Choose $\theta \in \ratQ^{\Irr G}$ with $\theta(\xxF)=0$ and let $\widetilde \theta$ be as in Definition \ref{thetatilde} with values \eqref{choicechikappa} of $\chi$ and $\kappa$. In this setting, do we have
$$ \theta(\xxF') \sgeq 0 \quad \Longleftrightarrow \quad \widetilde \theta(\xxF') \sgeq 0 \quad?$$
If not, are there additional assumptions on $\theta$ under which this equivalence holds?
\end{ques}
The fact that this is so hard to decide indicates that the passage from finite to infinite groups is a profound issue.

\subsection{Theory of Hilbert functions}

Regarding Proposition \ref{endlvieleHilbertfkt}, it seems that one has to study the properties of Hilbert functions extensively to answer Question \ref{questionequiv}. In particular, one should consider the following questions:
\begin{ques} Let $h\colon \Irr G \to \natN_0$ be the Hilbert function of some $G$--module such that $h$ is determined by the values $h(\rho)$ for $\rho$ in some finite subset $D_- \subset \Irr G$.
\begin{enumerate}
 \item Which kinds of functions are possible for $h$?
 \item Let $h'\colon \Irr G \to \natN_0$ be a function determined by the $h'(\rho)$ for $\rho$ in $D_-$ and $h'(\rho) \leq h(\rho)$ for every $\rho \in \Irr G$. If $h'$ occurs as a Hilbert function of a $G$--equivariant coherent subsheaf of a $(G,h)$--constellation, what are the possible values of $h'$?
\end{enumerate}
\end{ques}

\subsection{Resolution of singularities}

The original purpose of our construction of $\sMth$ was the search for resolutions of singularities, especially in the symplectic setting. Therefore, one would have to investigate the following:

\begin{ques}
 Is $\Mth$ or $\sMth$ smooth or does there exist a smooth connected component?
\end{ques}

\begin{ques}
Assume $h(\rho_0) = 1$. Is $\eta\colon \sMth \to X \red G$ projective?
\end{ques}

Further, we want to know:
\begin{ques}
Assume $h(\rho_0) = 1$. Is the map $\eta\colon\Mth \to X \red G$ or its restriction to a smooth connected component a resolution of singularities? If this is the case and if $X \red G$ is a symplectic variety, is $\eta$ even a symplectic resolution?
\end{ques}

Conversely, inspired by the situation for finite $G$ examined in \cite{CI:2004}, we can ask:
\begin{ques}\hfill\par
\begin{enumerate}
 \item Is every crepant resolution of singularities of $X \red G$ a component of some moduli space of $\theta$--stable $(G,h)$--constellations $\sMth$ for an appropriate choice of $\theta$?
 \item What is the relation between the spaces $\sMth$ for different choices of $\theta$? For example, is there a chamber structure in the space $\ratQ^{\Irr G}$ such that if $\theta'$ is in an adjacent wall of the chamber of $\theta$, then there exists a morphism $\sMth \to M_{\theta'}(X)$?
\end{enumerate}
\end{ques}

In particular, consider the action of $Sl_2$ on $(\kxC^2)^{\oplus 6}$ and its restriction to the zero fibre $\mu^{-1}(0)$ of the moment map. In \cite{Bec:2010} we determined the symplectic variety $\mu^{-1}(0) \red Sl_2$ to be a nilpotent orbit closure and we found two symplectic resolutions of singularities, namely the cotangent bundle $T^*\prP^3$ and its dual $(T^*\prP^3)^*$. In \cite{Bec:2011} we showed that these are dominated by a non--symplectic resolution, given by an invariant Hilbert scheme:
$$\begin{xy} \xymatrix{
 & \SlHilb{\mu^{-1}(0)} \ar@{->}[ld] \ar@{->}[dd] \ar@{->}[rd] & \\
T^*\prP^3 \ar@{->}[dr] & & (T^*\prP^3)^* \ar@{->}[dl]\\
 & \mu^{-1}(0) \red Sl_2 &
}\end{xy}$$
We have $\SlHilb{\mu^{-1}(0)} = \sMth$ for $\theta \in \ratQ^{\Irr G}$ such that $\theta_{\rho_0}$ is the only negative value. We would like to find out if the same is true for the symplectic resolutions:
\begin{ques}
Are the symplectic resolutions $T^*\prP^3$ and $(T^*\prP^3)^*$ also of the form $\sMth$ and if so, which is the correct choice for $\theta$?
\end{ques}

\end{document}